\documentclass[]{interact}

\usepackage{epstopdf}
\usepackage[caption=false]{subfig}
\usepackage[numbers,sort&compress]{natbib}
\bibpunct[, ]{[}{]}{,}{n}{,}{,}
\renewcommand\bibfont{\fontsize{10}{12}\selectfont}
\makeatletter
\def\NAT@def@citea{\def\@citea{\NAT@separator}}
\makeatother

\theoremstyle{plain}
\newtheorem{theorem}{Theorem}[section]
\newtheorem{lemma}[theorem]{Lemma}
\newtheorem{corollary}[theorem]{Corollary}
\newtheorem{proposition}[theorem]{Proposition}

\theoremstyle{definition}
\newtheorem{definition}[theorem]{Definition}
\newtheorem{example}[theorem]{Example}

\theoremstyle{remark}
\newtheorem{remark}{Remark}

\newcommand{\FF}{\mathcal F}
\newcommand{\LL}{\mathcal L}
\newcommand{\HH}{\mathcal{H}}
\newcommand{\supp}{{\rm supp\,} }
\newcommand{\be}{\beta}
\newcommand{\ga}{\gamma}
\newcommand{\la}{\lambda}

\newcommand{\al}{\alpha}
\newcommand{\eps}{\varepsilon}

\newcommand{\iv}{^{-1} }

\newcommand {\R} {\mathbb R}
\newcommand {\N} {\mathbb N}

\newcommand {\B} {\mathbb B}

\newcommand {\dom} {{\rm dom}\,}
\newcommand {\epi} {{\rm epi}\,}

\newcommand {\cl} {{\rm cl}\,}

\newcommand {\sd} {\partial}

\def\es{\emptyset}
\def\lsc{lower semicontinuous}

\def\LHS{left-hand side}
\def\RHS{right-hand side}

\newcommand{\norm}[1]{\left\Vert#1\right\Vert}
\newcommand{\abs}[1]{\left\vert#1\right\vert}

\newcounter{mycount}
\def\cnta{\setcounter{mycount}{\value{enumi}}}
\def\cntb{\setcounter{enumi}{\value{mycount}}}
\newcommand{\AK}[1]{\todo[inline]{AK {#1}}}


\renewcommand {\theenumi} {\rm(\roman{enumi})}

\begin{document}
	\title{Zero Duality Gap Conditions Via Abstract Convexity}
	\articletype{ARTICLE TEMPLATE}
	\author{
		\name{Hoa T. Bui\textsuperscript{a}\thanks{CONTACT Hoa T. Bui. Email: hoa.bui@curtin.edu.au} 
			and Regina S. Burachik\textsuperscript{b}  
			and Alexander Y. Kruger\textsuperscript{c} 
			and David T. Yost\textsuperscript{c}}
		\affil{\textsuperscript{a} ARC Centre for Transforming Maintenance through Data Science and School of Electrical Engineering, Computing, and Mathematical Sciences, Curtin University, Australia; \\
			\textsuperscript{b} Mathematics, STEM University of South Australia;
			\\
			\textsuperscript{c} School of Engineering, Information Technology  and Physical Science, Federation University Australia, POB 663, Ballarat, Vic, 3350, Australia}
		\vspace{5mm} {Dedicated to the memory of Alexander Rubinov, an outstanding mathematician, a great person, teacher and friend, on the occasion of his 80$^{th}$ birthday}}
	
	\maketitle
	
	\begin{abstract} 
		Using tools provided by the theory of {\em abstract convexity}, we extend conditions for zero duality gap to the context of nonconvex and nonsmooth optimization. Mimicking the classical setting, an {\em abstract convex} function is the upper envelope of a family of {\em abstract affine} functions (being conventional vertical translations of the {\em abstract linear} functions). We establish new conditions for zero duality gap under no topological assumptions on the space of abstract linear functions. In particular, we prove that the zero duality gap property can be fully characterized in terms of an inclusion involving (abstract) {\em $\varepsilon-$subdifferentials}. This result is new even for the classical convex setting. Endowing the space of abstract linear functions with the topology of pointwise convergence, we extend several fundamental facts of functional/convex analysis. This includes (i) the classical Banach--Alaoglu--Bourbaki theorem (ii) the subdifferential sum rule, and (iii) a  constraint qualification for zero duality gap which extends a fact established by Borwein, Burachik and Yao (2014) for the conventional convex case. As an application, we show with a specific example how our results can be exploited to show zero duality for a family of nonconvex, non-differentiable problems.
	\end{abstract}
	
	
	\begin{amscode}
		52A01; 47N10; 49J53; 49J27; 49J52; 90C30
	\end{amscode}
	\section{Introduction}\label{AbstractS1}
	
	The theory of \emph{abstract convexity}, also called \emph{convexity without linearity}, is a powerful tool that allows to extend many facts from classical convex analysis to more general frameworks. This theoretical framework (i) provides a fresh interpretation of existing notions, (ii) generates new links between previously disconnected tools from optimization, and (iii) allows one to analyze these tools within an elegant, unified setting.
	
	It has been the focus of active research for the last fifty years because of its many applications in functional analysis, approximation theory, and nonconvex analysis. Just like conventional convex analysis, the development of abstract convexity has been mainly motivated by applications to optimization; see \cite{Mar,Bru,Kay,RubUde01,YaoLi18,Iof001, DutMarRub04,DutMarRub04a,RubDza02,NedOzdRub07,Shv, DarMoh11,Rub01,EberMoh10,MohSam,DoaMoh09,DutMarRub08, RubWu07,BukRub07,SatMoh19-ref1,Wu2007-ref2,MarcoRub-ref3,Ivan99-ref4,RubGlover-ref5}.
	In particular, the work  \cite{Rubinov1999-ref6} uses the tools of abstract convexity to develop an implementable, bundle-type, global optimization algorithm. More global optimization methods based on abstract convexity tools can be found in \cite{MikAnd-ref7}. Applications of abstract convexity to set-valued analysis can be found in \cite{EberMoh10,BukRub07}. This fruitful theoretical framework also provides alternative interpretations of fundamental tools as local calmness (see \cite{article-ref11}), or relevant families of functions, such as lower semicontinuous radiant functions \cite{article-ref12} or topical and sub-topical functions \cite{Ref13}. As recent examples, the tools of abstract convexity are used in \cite{Syg16,Syg18,EwaSyg20,GorTyk19,GorTyk20} to derive criteria for global minima and maxima of nonsmooth functions, and stronger versions of Lagrangian duality and minimax theorems that are applicable to lower semicontinuous functions which are bounded below by a quadratic function. The recent work \cite{Ref9} derives necessary optimality conditions for nonsmooth nonconvex constrained optimization problems, and in \cite{Ref9} we find characterizations of global minima for certain families of (non-convex) problems \cite{GorTyk19,GorTyk20}.
	
	A deep study on abstract convexity can be found in the seminal book of Alexander Rubinov \cite{Rub},
	see also the monograph by Ivan Singer \cite{ISing}.
	The whole idea of abstract convexity originates in one of the fundamental facts of convex analysis: every lower semicontinuous convex function $f$ is the {\em upper envelope} of affine functions: for all $x$,
	\begin{equation}\label{Abstractclassical}
		f(x)=\sup\{h(x):\; h \text{ is an affine function},\; h \le f\}.
	\end{equation}
	Many key results in convex analysis are consequences of the two important aspects of \eqref{Abstractclassical}:
	\begin{enumerate}
		\item
		the ``supremum'' operation, and
		\item
		the set of functions over which this supremum is taken.
	\end{enumerate}
	Abstract convexity approach retains the ``supremum'' operation in aspect (i) but replaces linear, and consequently also affine, functions in (ii) by other families of functions called, respectively, {\em abstract linear} and {\em abstract affine}. Since aspect (i) of \eqref{Abstractclassical} is retained, global properties of convex analysis are preserved even when dealing with nonconvex models. This approach is sometimes called a ``non-affine global support function technique" (see, for example, \cite{Rub,Kay,BukRub07}).
	
	The key tools from convex analysis, such as {\it subdifferentials}, $\varepsilon$-{\it subdifferentials} and {\it Fen\-chel--Mo\-reau conjugates} have their ``abstract'' counterparts, constructed by using abstract linear functions. For instance, the  abstract subdifferential of an abstract convex function $f$ at a point $x$ collects all the supporting abstract linear functions which are minorants of $f$ (i.e., their graphs stay below the graph of $f$), and coincide with $f$ at $x$.
	It extends the concept of convex subdifferential and provides a valuable tool for studying certain nonconvex optimization problems; see \cite{Rub,Kay,BukRub07,RubWu07}.
	
	In \cite{Bor}, zero duality gap is shown to be equivalent to (a) certain properties involving $\varepsilon$-sub\-differentials and (b) certain facts involving conjugate functions. One of the aims of the present work is to extend these results to the context of abstract convexity. Additionally, we supplement the sum rule for abstract subdifferentials, improving the corresponding result in \cite{Kay}. To the best of our knowledge, the only attempt to consider explicitly the pointwise convergence topology on the space of abstract linear functions to deduce calculus rules for subdifferentials was made in \cite{Bru}. Here we exploit this idea further. In particular, we extend the fundamental Banach--Alaoglu--Bourbaki theorem on the weak$^*$ closedness of the dual unit ball to the general space of abstract linear functions. To apply the new theory, we show how zero duality gap can be established for a family of nonconvex and nonsmooth problems. 
	
	The structure of the paper is as follows. Section \ref{AbstractS2} recalls some preliminary definitions and facts used throughout the paper. We briefly introduce and study the space of abstract linear functions, and abstract convexity notions. Some results in this section are new.
	In Section \ref{AbstractS3}, we provide properties which ensure the zero duality gap. We impose no topological assumptions on the primal space nor on the space of linear functions. A comparison with its forerunner, \cite[Theorem 3.2]{Bor} for convex programming, is established. A necessary and sufficient characterization of zero duality gap is provided, which is new even in the standard convex setting. In Section \ref{S4}, we equip the abstract linear function space with the pointwise convergence topology to extend some classical convex subdifferentials' calculus in the framework of abstract convexity. Here, we assume that the epigraphs of conjugate functions admit a certain additivity property (see \eqref{Abstractcon}). This condition holds for conventional lower semicontinuous convex functions. The Banach--Alaoglu--Bourbaki theorem is extended to a general function space.
	Then the conditions for the zero duality gap are exposed fully.
	In Section \ref{AbstractExample_AC}, we construct a nontrivial example for which our analysis and the zero duality gap property apply. Finally, Section \ref{sec-oq} lists some open questions.
	
	\section{Preliminaries on Abstract Convexity}\label{AbstractS2}
	
	In what follows, $\R$ and $\N$ stand for the sets of all real numbers and all positive integers, respectively.
	We use the notation $\R_{+\infty}:=\R\cup\{+\infty\}$ and $\R_{\pm\infty}:=\R\cup\{\pm\infty\}$.
	Throughout, $X$ is a nonempty set. Unless explicitly stated otherwise, we do not assume any algebraic or topological structure on $X$. Given a function $f:X\to\R_{\pm\infty}$, its domain and epigraph are the sets $\dom f:=\{x\in X:\, f(x)<+\infty\}$ 
	
	The sum of functions from $X$ to $\R_{+\infty}$ is defined in the usual way: $(f_1+f_2)(x):=f_1(x)+f_2(x)$ for all $x\in X$, and {we} write $f_1\le f_2$ if $f_1(x)\le f_2(x)$ for all $x\in X$. {For $A,B$ nonempty sets, and a point-to-set mapping $T:A\rightrightarrows B$, we define ${\rm Dom\,} T:=\{a\in A\::\: Ta \neq\emptyset\}$. Given a subset $D\subset A$ the {\em direct image of $D$ by $T$} is the set $T(D):=\cup_{d\in D} Td$.}
	
	When talking about convex functions or convex sets,
	we mean the
	conventional
	convexity.
	We recall next the standard strict convex separation theorem; cf., e.g., \cite[Theorem~1.7]{Bre}.
	\begin{lemma} \label{AbstractS.G.F}
		Let $A$ and $B$ be nonempty convex subsets of a normed vector space such that $A$ is closed and $B$ is compact. If $A\cap B=\emptyset$, then there is a nonzero linear continuous functional $x^*$ such that $\sup_{x\in A} \langle x^*, x\rangle < \inf_{y\in B} \langle x^*,y\rangle.$
	\end{lemma}
	As mentioned above, linear functions and their vertical shifts (affine functions) are at the core of convex analysis. They
	play
	a crucial role in the definitions of subdifferentials and conjugate functions.
	In the next subsection, we define the
	space
	$\LL$ of abstract linear functions
	and show
	that some classical features of subdifferentials and conjugate functions
	remain
	true for this more general set of linear functions. Most definitions and notation arising from abstract convexity theory are standard and taken from \cite{Rub}.

	\subsection{Space of Abstract Linear functions}
	
	Throughout,
	in addition to the given set $X$,
	we assume that $\LL$
	is a
	given
	space of abstract linear functions. We may alternatively call $\LL$ the {\em abstract linear space.}
	
	\begin{definition}\label{AbstractALS}
		A {\em space $\LL$ of abstract linear functions}
		is  a
		family of functions $l:X\to\R$
		satisfying the following properties:
		
		
		\begin{enumerate}
			\item
			$\LL$ is closed with respect to the addition
			operation,
			i.e. $l_1,l_2\in \LL \Longrightarrow l_1+l_2\in \LL$;
			\item
			for every $l\in \LL$ and $m\in \N$, there exist $l_1,\ldots,l_m\in \LL$ such that
			$l=\sum_{i=1}^m l_i$.
		\end{enumerate}
	\end{definition}
	
	\begin{remark}
		\begin{enumerate}
			\item
			Assume also that  $\LL$ verifies
			part (i) of Definition \ref{AbstractALS}, and that it also has a neutral element of the addition. Namely, there exists $0\in \LL$ such that $0+l=l$ for all $l\in\LL$. Then $\LL$ automatically verifies
			part (ii).
			
			\item
			Note that $\LL$ is not in necessarily a vector space because scalar multiplication is not required.
		\end{enumerate}
	\end{remark}
	
	
	
	\begin{definition}\label{AbstractD3.1}
		\begin{enumerate}
			\item
			Let $X$ be
			equipped with an addition operation,
			and $m\in\N$.
			The {\it infimal convolution}
			of
			functions $\psi_1,\ldots,\psi_m: X\to \R_{+\infty}$ is the function 
			\begin{equation*}
				\psi_1\square \ldots\square \psi_m(x):= \inf_{x_1+\ldots+x_m=x}\left\{\psi_1(x_1)+\ldots+\psi_m(x_m) \right\},\quad x\in X,
			\end{equation*}
			with the convention that infimum over
			the
			empty set is $+\infty$. {We say that the infimal convolution is {\em exact} when the infimum in the expression above is attained.}
			\item
			The {\it (Fenchel) conjugate} of a function $f:X\to \R_{+\infty}$ is 
			the function 
			\begin{equation*}
				f^*(l):=\sup_{x\in X}\{l(x)-f(x)\}, \quad l\in\LL.
			\end{equation*}
			{Similarly, the {\it bi-conjugate} of $f$ is 
				the function $f^{**}: X\to \R_{+\infty}$, defined by
				\begin{equation*}
					f^{**}(x):=\sup_{l\in \LL}\{l(x)-f^*(l)\}, \quad x\in X.
				\end{equation*}
			}
			\item
			Given
			a number $\eps \ge 0$,
			the {\it $\eps-$subdifferential}
			of a function $f:X\to \R_{+\infty}$
			at a point $x\in \dom f$
			is the point-to-set mapping defined as
			\begin{equation*}
				\partial_{\eps} f(x):=\{l\in \LL: f(y)-f(x)-(l(y)-l(x))+\eps \ge 0 \text{ for all } y\in X\}.
			\end{equation*}
			If $x\notin \dom f$, then $\partial_\eps f(x)=\emptyset$.
			If $\eps=0$, we say simply `subdifferential' and write $\sd f(x)$.
		\end{enumerate}
	\end{definition}

	\begin{remark}
		Given a function $f:X\to \R_{+\infty}$, it follows from the definition of $\eps-$subdifferential that
		\begin{equation}\label{domf2}
			\displaystyle\bigcap_{\eps>0} {\rm Dom\,}\partial_{\eps} f \subset \dom f.
		\end{equation}
	\end{remark}
	
	
	The next proposition summarizes the key properties of the concepts {given} in Definition~\ref{AbstractD3.1}.
	In the conventional convex setting, these properties are well known.
	
	\begin{proposition}\label{AbstractP3.2}
		Suppose $f:X\to \R_{+\infty}$.
		The following
		assertions
		hold.
		\begin{enumerate}
			\item
			For all $x\in\dom f$
			and $\eps \ge 0$,
			$l\in \partial_{\eps}f(x)$ if and only if
			$f^*(l)+f(x)\le l(x)+\eps$.
			\item
			\[
			\dom f^*= \bigcap_{\eps >0} \partial_\eps f(X).
			\]
			\item
			For all $x\in X$
			and $\eps \ge 0$,
			\[
			\bigcap_{\eta>0}\partial_{\eps+\eta} f(x)=\partial_{\eps} f(x).
			\]
			{ In particular, $\partial f(x)= \bigcap_{\eps>0}\partial_{\eps} f(x).$}
			\cnta
		\end{enumerate}
		\noindent Suppose $f_1,\ldots,f_m:X\to \R_{+\infty}$ $(m\ge 2)$.
		\begin{enumerate}
			\cntb
			\item
			$\left(\sum_{i=1}^m f_i\right)^*
			\le
			f_1^*\square\ldots\square f_m^*$.
			
			\item
			%
			For all $x\in X$ and $\eps \ge 0$,
			\begin{equation}
				\label{AbstractP3.2.v}
				\bigcup_{\substack{
						\eps_i\ge 0
						,\;i=1,\ldots,n
						\\
						\sum_{i=1}^m\eps_i=\eps
				}}\;
				\sum_{i=1}^m \partial_{\eps_i}f_i(x)
				\subset
				\partial_{\eps}\left(\sum_{i=1}^mf_i\right)(x),
			\end{equation}
			and consequently,
			\begin{equation}
				\label{AbstractP3.2.v2}
				\bigcap_{\eta >0} \;
				\bigcup_{\substack{
						\eps_i\ge 0
						,\;i=1,\ldots,n
						\\
						\sum_{i=1}^m\eps_i=\eps+\eta
				}}\;
				\sum_{i=1}^m \partial_{\eps_i}f_i(x)
				\subset
				\partial_{\eps}\left(\sum_{i=1}^mf_i\right)(x).
			\end{equation}
		\end{enumerate}
	\end{proposition}
	
	\begin{proof}
		\begin{enumerate}
			\item
			See \cite[Proposition 7.10]{Rub}.
			\item
			Recall that $\partial_\eps f(X)=\bigcup_{x\in X} \partial_\eps f(x)$.
			Let $l\in \dom f^*$ and fix $\eps >0$. { We need to find $x\in X$ such that $l\in \partial_{\eps} f(x)$.}
			By Definition~\ref{AbstractD3.1}(ii), there is an $x\in X$
			such that $l(x)-f(x)> f^*(l)-\eps$.
			Hence,
			$f^*(l)+f(x)< l(x)+\eps$, and it follows from (i) that $l\in \partial_{\eps} f(x)$. { Since $\eps >0$ is arbitrary,
				this proves that $l\in \cap_{\eps>0} \cup_{x\in X} \partial_{\eps} f(x)$.}
			Conversely, let $\eps >0$ and $l\in \partial_{\eps} f(x)$ for some $x\in X$.
			By (i), $f^*(l)\le l(x)-f(x)+\eps<+\infty$, i.e. $l\in \dom f^*$. { Hence, the 
				{equality 
					holds}.}
			\item { The last statement in (iii) follows directly from the the first one. Hence, we proceed to prove the first statement.}
			Let $x\in X$ and $\eps \ge 0$.
			It is clear from Definition~\ref{AbstractD3.1}(iii) that
			$\partial_{\eps} f(x)\subset\bigcap_{\eta>0}\partial_{\eps+\eta} f(x)$.
			{To prove}
			the opposite inclusion, let $l\in \bigcap_{\eta>0}\partial_{\eps+\eta}f(x)$.
			Then, by (i), $f^*(l)+f(x)\le l(x)+\eps+\eta$
			for all $\eta >0$,
			and consequently, $f^*(l)+f(x)\le l(x)+\eps$.
			Using (i)  again, we
			conclude
			that $l\in \partial_{\eps} f(x)$.
			\item
			Let $l\in \LL$.
			Take any additive decomposition of $l$,
			i.e. a finite collection $l_1,\ldots,l_m\in \LL$ such that
			$l_1+\ldots+l_m=l$.
			We consider two cases.
			
			{\bf Case 1.}
			Suppose $l\notin \dom \left(\sum_{i=1}^m f_i\right)^*$.
			By Definition~\ref{AbstractD3.1}(ii),
			\begin{align*}
			\sum_{i=1}^m f_i^*(l_i)
			&=
			\sum_{i=1}^m\sup_{x\in X}(l_i(x)-f_i(x))\\
			&\ge
			\sup_{x\in X}\left\{\sum_{i=1}^m(l_i(x)-f_i(x)) \right\}=\left(\sum_{i=1}^m f_i\right)^*(l)=+\infty.
			\end{align*}
			{Since this is true for every additive decomposition of $l$, } by Definition~\ref{AbstractD3.1}(i),
			$f_1^*\square \ldots\square f_m^*(l)=+\infty.$  {Note that in this case the infimal convolution is trivially exact.}

			{\bf Case 2.}
			Suppose $l\in \dom \left(\sum_{i=1}^m f_i\right)^*$.
			By Definition~\ref{AbstractD3.1}(ii),
			for any $\eps >0$,
			there is an $x\in X$ such that
			\begin{align*}
				\left(\sum_{i=1}^m f_i\right)^*(l)&\le l(x)-\sum_{i=1}^m f_i(x)+\eps= \sum_{i=1}^m(l_i(x)-f_i(x))+\eps\le \sum_{i=1}^mf_i^*(l_i)+\eps,
			\end{align*}
			{where the additive decomposition of $l$ was used in the equality. Hence, in view of Definition~\ref{AbstractD3.1}(i), and taking infimum over all possible additive decompositions and 
				{all}
				$\eps>0$, we obtain}
			$$
			\left(\sum_{i=1}^m f_i\right)^*(l)
			\le \inf_{l_1+\ldots+l_m=l,\,\eps>0} \left\{\sum_{i=1}^mf_i^*(l_i)+\eps\right\}
			=f_1^*\square \ldots \square f_m^* (l),
			$$
			establishing 
			assertion
				(iv).
			%
			\item
			{Note first that it is enough to prove \eqref{AbstractP3.2.v}. 
				Indeed, \eqref{AbstractP3.2.v2} follows directly from taking $\cap_{\eta>0}$ on both sides of inclusion \eqref{AbstractP3.2.v} with $\eps +\eta$ in place of $\eps$, combined with part (iii) of this proposition. Hence, we proceed to establish \eqref{AbstractP3.2.v}.}
			Let $x\in X$, $\eps \ge 0$,
			$\eps_i\ge 0$, $l_i\in \partial_{\eps_i}f_i(x)$ $(i=1,\ldots,m)$,
			$\sum_{i=1}^m\eps_i=\eps$ and $\sum_{i=1}^m l_i=l$.
			By (i), $f^*_i(l_i)+f_i(x)\le l_i(x)+\eps_i$
			$(i=1,\ldots,m)$. {Using also the definition of infimal convolution, we obtain}
			\[
			\begin{array}{rcl}
				f_1^*\square \ldots \square f_m^* (l) +  \sum_{i=1}^m f_i(x) &\le& \sum_{i=1}^m f^*_i(l_i)+\sum_{i=1}^m f_i(x)
				\\
				&&\\
				&=&\sum_{i=1}^m f^*_i(l_i)+ f_i(x) \le \sum_{i=1}^m l_i +\eps_i =  l(x)+\eps.
			\end{array}
			\]
			{Using now (iv) in the left hand side of the expression above, we deduce}
			\begin{equation*}
				\left(\sum_{i=1}^m f_i\right)^*(l)+\left(\sum_{i=1}^m f_i\right)(x)
				\le l(x)+\eps,
			\end{equation*}
			which, in view of (i), is equivalent to $l\in\sd_{\eps}(\sum_{i=1}^m f_i)(x)$.
			{This proves \eqref{AbstractP3.2.v}, as wanted.}
		\end{enumerate}
	\end{proof}
	
	\subsection{Abstract Convex Functions and Abstract Convex Sets}
	
	Given the space $\LL$ of abstract linear functions on $X$, we now introduce another fundamental {concept} of abstract convex analysis: the space $\HH$ of abstract affine functions.
	
	\begin{definition}\label{AbstractD4.1}
		The space $\HH$ of {\em abstract affine functions} is defined as
		$\HH:=\{l+c: l\in \LL,\; c\in \R\}.$
	\end{definition}
	
	\begin{remark}\label{R7}
		\begin{enumerate}
			\item
			Obviously $\LL\subset\HH$.
			\item
			The space of abstract affine functions
			can be defined independently
			of
			the space of abstract linear functions
			as
			a
			set of functions which is closed with respect to
			vertical shifts.
			\item
			The space of abstract affine functions satisfies the conditions in Definition~\ref{AbstractALS} and as such can also be considered as a space of abstract linear functions.
			As a consequence, many facts formulated in the paper for the space $\LL$ are also valid for the space $\HH$.
		\end{enumerate}
	\end{remark}
	
	
	We can now extend the classical notions of convex functions and convex sets to
	the
	abstract convexity framework;
	cf. \cite[Definition 7.2]{Rub}.
	
	\begin{definition}\label{AbstractD4.2}
		
		\begin{enumerate}
			\item
			Given a function $f:X\to \R_{+\infty}$,
			its {\em support set} is given by
			\begin{equation*}
				\supp f:=\{h\in \HH: h(x)\le f(x) \text{ for all } x\in X \}.
			\end{equation*}
			{Namely, $\supp f$ is the set of all minorants of $f$ which belong to $\HH$.}
			\item
			A function $f:X\to \R_{+\infty}$ is
			\emph{$\LL-$convex} if
			there is a subset $L\subset \LL$ such that
			$f(x)=\sup_{l\in L}l(x)$ for all $x\in X$.
			\item
			A function $f:X\to \R_{+\infty}$ is
			{\it $\HH-$convex}
			if there is a
			subset $H\subset \HH$ such that
			$f(x)=\sup_{h\in H}h(x)$ for all $x\in X$.
			
			\item
			A subset $C\subset \LL$ is
			\emph{$\LL-$convex}
			if for any
			$l_0\notin C$,
			%
			there is an $x\in X$ such that
			$
			l_0(x)>\sup_{l\in C}l(x).
			$
			\item
			A subset $C\subset \HH$ is $\HH-$convex if for any $h_0\notin C$ there is an $x\in X$ such that
			$
			h_0(x)>\sup_{h\in C} h(x).
			$
		\end{enumerate}
	\end{definition}

	
	\begin{remark}\label{Abstract-R}
		\begin{enumerate}
			\item {Since we restrict our analysis to functions from $X$ to $\R_{+\infty}$ (i.e., $f>-\infty$), the sets $L$ and $H$ in 
				{parts (ii) and (iii) of Definition \ref{AbstractD4.2}}
				cannot be empty.} 
			\item
			The definition of $\LL-$convex set in Definition~\ref{AbstractD4.2}(iv)
			differs
			from
			\cite[Definition 1.4]{Rub}
			(as well as the one used in \cite{Kay})
			although is equivalent to it; see
			Proposition~\ref{AbstractP4.2.2}(i).
			The two versions of the definition reflect
			different ideas of convexity. Definition~\ref{AbstractD4.2}(iv)
			exploits the {\it separation} property, whereas the rationale in
			\cite[Definition 1.4]{Rub}
			is based on the notion of convex combination.
			\item
			{The definition implies that for any $f:X\to \R_{+\infty}$,
				its support set $\supp f$ is $\HH-$convex. Indeed, take any $\hat{h}\notin \supp f$. By definition of support set there exists $\hat{x}\in X$ such that $\hat{h}(\hat{x}) > f(\hat{x})\ge \sup_{h\in \supp f} h(\hat{x})$, which means that $\supp f$ is $\HH-$convex by Definition~\ref{AbstractD4.2}(iii).}  
			\item
			{Property (iii) also goes in the opposite way: If $C$ is any $\HH-$convex set and $f_C(\cdot):=\sup_{l\in C} l(\cdot)$, then $C=\supp f$. Indeed, the definition of $f_C$ readily gives $C\subset \supp f$. For the opposite inclusion, assume there exists $\hat{h}\in \supp f$ such that $\hat{h}\notin C$. Since $C$ is $\HH-$convex there exists $\hat{x}\in X$ such that $\hat{h}(\hat{x}) > \sup_{h\in C} h(\hat{x})=f_C(\hat{x})\ge \hat{h}(\hat{x})$, where the last inequality holds because $\hat{h}\in \supp f$. This contradiction implies that 
				$C=\supp f$.}
			\item {Because every $l\in \LL$ is finite valued, so is every $h\in \HH$. Hence,  $h\in \supp f$ if and only if $h(x)\le f(x)$ for every $x\in \dom f$.}
		\end{enumerate}
	\end{remark}
	
	The
	next proposition collects some properties of $\HH-$convex functions.
	
	\begin{proposition}\label{AbstractP4.1}
		Suppose
		$f:X\to \R_{+\infty}$. The following assertions hold.
		\begin{enumerate}
			\item
			For all $x\in X$, 
			\begin{equation}\label{AbstractP4.1.1}
				\sup_{h\in \supp f}h(x)\le f(x).
			\end{equation}
			Equality holds in \eqref{AbstractP4.1.1} for all $x\in X$ if and only if $f$ is $\HH-$convex.
			\item
			$(l,r)\in\epi f^*$ if and only if $l-r\in \supp f$.
			\item  $f$ is $\HH-$convex if and only if
			\begin{equation}\label{domf22}
				\dom f=\bigcap_{\eps>0} {\rm Dom\,} \partial_{\eps} f.
			\end{equation}
			Consequently,  if $f$ is $\HH-$convex, then
			$\partial_{\eps}f(x)\neq \emptyset$
			for all $x\in \dom f$ and $\eps >0$.
			{Vice versa,}
			$f$ is $\HH-$convex if
			$\partial_{\eps}f(x)\neq \emptyset$
			for all {$x\in \dom f$} and $\eps >0$.
			\item
			If $x\in \dom f$ and $(l,r)\in \epi f^*$, then $\eps:=r+f(x)-l(x)\ge0$ and
			$l\in \partial_{\eps}f(x)$.
			\item (Fenchel--Moreau) For all $x\in X$, 
			\begin{equation}
				\label{AbstractP4.1.2}
				f^{**}(x)\le f(x).
			\end{equation}
			Equality holds in \eqref{AbstractP4.1.2} for all $x\in X$ if and only if $f$ is $\HH-$convex.
		\end{enumerate}
	\end{proposition}
	
	\begin{proof}
		\begin{enumerate}
			\item
			Inequality \eqref{AbstractP4.1.1} holds trivially by the definition of $\supp f$.
			If
			equality holds in \eqref{AbstractP4.1.1} for all $x\in X$, then $f$ is $\HH-$convex by Definition~\ref{AbstractD4.2}(iii), with { the subset $H:=\supp f$}.
			Conversely, assume that $f$ is $\HH-$convex.
			By
			Definition~\ref{AbstractD4.2}(iii),
			there exists a set $H\subset \HH$ such that
			$f(x)=\sup_{h\in H}h(x)$
			for all $x\in X$,
			and consequently,
			$H\subset \supp f$.
			Hence,
			$$f(x)=\sup_{h\in H}h(x)\le \sup_{h\in \supp f}h(x)\le f(x).$$
			This proves that equality holds in \eqref{AbstractP4.1.1}. 
			\item
			{The proof of this fact can be found in} \cite[page 444, equation (1)]{Kay}. {For the reader's convenience, we sketch here the proof: $(l,r)\in\epi f^*$ if and only if $r\ge f^*(l)\ge l(x)-f(x)$ for every $x\in X$. The latter can be re-arranged as $l(x)-r\le f(x)$  for every $x\in X$. Equivalently,
				$l(\cdot)-r\in \supp f$.}
			
			\item
			{Note first that, by definition, the right hand side of \eqref{domf22} is always a subset of the left hand side. Hence, we only need to show the opposite inclusion. Assume that the function $f$ is $\HH-$convex and take $x\in \dom f$. Fix $\eps>0$ arbitrary.  By  (i), we have $f(x)-\eps <\sup_{h\in \supp f}h(x)=f(x)$.
				Then, there
				exists
				$l\in \LL$ and $c\in \R$ such that
				$$
				l(y)+c\le f(y) \text{ for all } y\in X
				\quad\text{ and} \quad
				l(x)+c+\eps> f(x).
				$$
				Consequently,
				$f(y)-f(x)> l(y)-l(x)-\eps$ for all $y\in X$,
				i.e.
				$l\in \partial_{\eps}f(x)$. Hence, $x\in {\rm Dom\,} \partial_{\eps}$. Since $\eps>0$ is arbitrary, we deduce that
				$x\in \bigcap_{\eps>0} {\rm Dom\,} \partial_{\eps} f$. Altogether, we obtained that $\dom f \subset \bigcap_{\eps>0} {\rm Dom\,} \partial_{\eps} f$
			}.
			Conversely,
			{let $x\in \bigcap_{\eps>0} {\rm Dom\,} \partial_{\eps} f$. Hence, for every $\eps>0$ there exists $l_{\eps}\in \partial_{\eps}f(x)$. Altogether, $f(y)-f(x)\ge l_{\eps}(y)-l_{\eps}(x)-\eps$ for all $y\in X$.
				Set $h_{\eps}(y):=l_{\eps}(y)-l_{\eps}(x)+f(x)-\eps$ ($y\in X$).
				Then $h_{\eps}\in\supp f$ and therefore $f(x)=h_{\eps}(x)+\eps\le\sup_{h\in \supp f}h(x)+\eps$.
				It follows that $f(x)\le\sup_{h\in \supp f}h(x)$ and, thanks to (i), $f(x)=\sup_{h\in \supp f}h(x)$.
				It now follows from (i) that
				$f$ is $\HH-$convex. The last two 
				{assertions}
				in (iii) are direct consequences of \eqref{domf22}.
			}
			\item
			Let $x\in \dom f$ and $(l,r)\in \epi f^*$.
			Then $r\ge f^*(l)\ge l(x)-f(x)$. Define $\eps:=r+f(x)-l(x)\ge0$.
			Hence,
			$f^*(l)+f(x)\le r+f(x)= l(x)+\eps$. Now Proposition~\ref{AbstractP3.2}(i) yields $l\in \partial_{\eps}f(x)$.
			\item
			See \cite[Theorem 7.1]{Rub}.
		\end{enumerate}
	\end{proof}
	
	\begin{remark}
		Inequality \eqref{AbstractP4.1.1} holds
		even when $\supp f=\emptyset$,
		as in this case
		{$f(x)=\sup_{h\in \supp f}h(x)=-\infty$.
			Since we assume that 
			{$f$ never equals}
			$-\infty$, this case is not considered. 
			Consequently, if $f:X\to \R_{+\infty}$ is $\HH-$convex,
			then necessarily
			$\supp f\not=\emptyset$.}
	\end{remark}
	
	The next proposition provides some properties of $\LL-$
	and $\HH-$convex sets
	used in the subsequent sections.
	
	\begin{proposition}\label{AbstractP4.2.2}
		\begin{enumerate}
			\item
			A
			nonempty
			subset $C\subset \LL$ is $\LL-$convex if and only if there is an $\LL-$convex function $f:X\to\R_{+\infty}$ such that
			$
			C=\supp f.
			$
			In this case, $\supp f\subset\LL$, i.e.
			$$\supp f= \{l\in \LL: l(x)\le f(x)\text{ for all } x\in X\},$$
			and consequently,
			$C=
			\{l\in \LL: f^*(l)\le 0\}$.
			
			\item
			A nonempty subset $C\subset \HH$ is $\HH-$convex if and only if there is an $\HH-$convex function $f:X\to\R_{+\infty}$ such that
			$C=\supp f$.
			\item
			Suppose
			$\LL$
			is a vector space
			(i.e. {it is closed with respect to addition and multiplication by scalars}).
			If a subset of $\LL$ is $\LL-$convex (a subset of $\HH$ is $\HH-$convex), then it is convex in the conventional sense.
			
			
		\end{enumerate}
	\end{proposition}
	
	\begin{proof}
		\begin{enumerate}
			\item See  \cite[Lemma 1.1]{Rub}.
			%
			The last representation follows in view of Definition~\ref{AbstractD3.1}(ii) of the conjugate function.
			\item
			is a consequence of (i) in view of Remark~\ref{R7}(iii).
			\item
			Suppose $C\subset\LL$ is $\LL-$convex and $f:X\to\R_{+\infty}$ is an $\LL-$convex function such that
			$C=\supp f$ (cf. (i)).
			Take $l_1,l_2\in C$ and $\alpha\in [0,1]$.
			Hence, for all $x\in X$, we have
			$l_1(x)\le f(x)$ and $l_2(x)\le f(x)$,
			and consequently,
			$\alpha l_1(x)+(1-\alpha)l_2(x)\le f(x)$,
			i.e.
			$\alpha l_1+(1-\alpha) l_2\in \supp f=C$.
			Since $\LL$ is a vector space, $\HH$ is a vector space too.
			The case of an $\HH-$convex set follows in view of Remark~\ref{R7}(iii).
		\end{enumerate}
	\end{proof}
	
	\begin{remark}
		{
			The converse of Proposition~\ref{AbstractP4.2.2}(iii) is not true:
			not all convex sets are $\LL-$convex;
			cf.
			the characterisations of $\LL-$convex sets
			in Proposition \ref{AbstractP4.E0}(ii), Section~\ref{AbstractExample_AC}.
		}
	\end{remark}

	\section{Conditions for Zero Duality Gap}\label{AbstractS3}
	
	We
	consider the minimization problem:
	\begin{align}
		\label{AbstractP}\tag{P}
		\inf\sum_{i=1}^m f_i(x)\quad
		\textnormal{s.t.}\quad x\in X,
	\end{align}
	where $X$ is a general nonempty set, and 
	$f_1,\ldots,f_m:X\to \R_{+\infty}$ ($m\ge 2$)
	are arbitrary functions.
	The corresponding dual problem has the following form:
	\begin{align}
		\label{AbstractD}\tag{D}
		\sup\sum_{i=1}^m(-f^*_i(l_i))\quad
		\textnormal{s.t.}\quad l_1,\ldots,l_m\in \LL,\;\;
		l_1+\ldots+l_m=0.
	\end{align}
	
	Denote by $v(P)$ and $v(D)$, the optimal values of \eqref{AbstractP} and \eqref{AbstractD}, respectively.
	In general, we have the inequality $v(P)\ge v(D)$.
	We say that a {\em zero duality gap} holds for problems \eqref{AbstractP} and \eqref{AbstractD} if $v(P)=v(D)$. We refer the reader to \cite{Bor} for comments and further explanation on the zero duality gap. The following
	relations in terms of the conjugate functions and infimal convolution are
	direct consequences of Definition~\ref{AbstractD3.1}(ii) and (i), respectively:
	\begin{align}
		\label{AbstractP1}
		v(P)=&\inf_{x\in X}\sum_{i=1}^m f_i(x)=-\left(\sum_{i=1}^m f_i\right)^*(0),\\
		\label{AbstractD1}
		v(D)=&\sup_{l_1+\ldots+l_m=0} \sum_{i=1}^m(-f^*_i(l_i)) =-(f_1^*\square \ldots\square f_m^*)(0).
	\end{align}
	Thus, zero duality gap is equivalent to
	\begin{equation}
		\label{AbstractEx-Zero}
		\left(\sum_{i=1} f_i\right)^*(0)=(f_1^*\square \ldots\square f_m^*)(0).
	\end{equation}
	
	The next theorem
	extends  \cite[Theorem 3.2]{Bor} to our general framework. It
	contains subdifferential characterizations of a stronger
	condition
	\if{
		$$
		\left(\sum_{i=1}^m f_i\right)^*=f_1^*\square \ldots\square f_m^*,
		$$
	}\fi
	which clearly
	ensures
	\eqref{AbstractEx-Zero}.
	The proof below refines the core arguments in the proof of \cite[Theorem 3.2]{Bor}.
	
	\begin{theorem}\label{AbstractT3.3}
		Let $f_1,\ldots,f_m,: X\to \R_{+\infty}$ $(m\ge 2)$ be such that $\bigcap_{i=1}^m\dom f_i\neq \emptyset$.  The following
		properties
		are equivalent:
		\begin{enumerate}
			\item
			for all $x\in X$, $\eps \ge 0$ and $K>1$, it holds
			\begin{gather}
				\label{AbstractT3.3.1}
				\partial_\eps \left(\sum_{i=1}^{m}f_i\right)(x)\subset \sum_{i=1}^m \partial_{K\eps}f_i(x);
			\end{gather}
			
			\item
			there is a $K>0$ such that inclusion \eqref{AbstractT3.3.1} holds for all $x\in X$ and $\eps\ge 0$;
			\item
			$\left(\sum_{i=1}^m f_i\right)^*=f_1^*\square\ldots\square f_m^*$.
			
			\item
			for all $x\in X$ and all $\eps \ge 0$,
			\begin{gather}
				\label{AbstractT3.3.2}
				\partial_\eps \left(\sum_{i=1}^{m}f_i\right)(x)= \bigcap_{\eta>0}\bigcup_{\substack{\eps_i\ge 0
						,\;i=1,\ldots,m
						\\
						\sum_{i=1}^m\eps_i=\eps+\eta}}\sum_{i=1}^m \partial_{\eps_i}f_i(x);
			\end{gather}

		\end{enumerate}
	\end{theorem}
	
	\begin{proof}
		\begin{enumerate}
			\item[] The implication (i) $\Rightarrow$ (ii) is obvious.
			
			\item[] (ii) $\Rightarrow$ (iii).
			Thanks to Proposition~\ref{AbstractP3.2}(iv), we only need to show that,
			for all ${l\in \LL}$,
			\begin{gather}
				\label{T13P1}
				(f_i^*\square \ldots\square f_m^*)(l)\le \left(\sum_{i=1}^m f_i\right)^*(l).
			\end{gather}
			If $l\notin \dom \left(\sum_{i=1}^mf_i\right)^*$, the inequality holds trivially.
			Let $l\in \dom \left(\sum_{i=1}^mf_i\right)^*$. Thanks to Proposition~\ref{AbstractP3.2}(ii), {for every $\eta> 0$ there is an $x_{\eta}\in X$ such that $l\in \partial_{\eta}\left(\sum_{i=1}^m f_i \right)(x_{\eta})$.
				By (ii), there are $l_i\in \partial_{K\eta} f_i(x_{\eta})$ ($i=1,\ldots,m$) such that $l=\sum_{i=1}^ml_i$.
				Thanks to Proposition~\ref{AbstractP3.2}(i),
				$$
				f^*_i(l_i)+f_i(x_{\eta})\le l_i(x_{\eta})+K\eta, \quad  i=1,\ldots,m.
				$$
				Using the definitions, we can write
				\begin{align*}
					(f_i^*\square \ldots\square f_m^*)(l)&\le \sum_{i=1}^m f_i^*(l_i)\le -\sum_{i=1}^m f_i(x_{\eta}) + \sum_{i=1}^m l_i(x_{\eta})+mK\eta\\
					&=-\left(\sum_{i=1}^m f_i\right)(x_{\eta})+l(x)+mK\eta \le \left(\sum_{i=1}^m f_i\right)^*(l)+mK\eta.
				\end{align*}
				Passing to the limit as $\eta\downarrow0$, we arrive at \eqref{T13P1}.
			}
			\item[] (iii) $\Rightarrow$ (iv).
			{Note first that when $x\not\in \bigcap_{i=1}^m\dom f_i$, both sides of \eqref{AbstractT3.3.2} are empty, so it is enough to prove the statement for $x\in \bigcap_{i=1}^m\dom f_i$. }Let $x\in \bigcap_{i=1}^m\dom f_i$ and $\eps \ge0$.
			Thanks to Proposition~\ref{AbstractP3.2}(v), we only need to show that
			\begin{gather}
				\label{T13P2}
				\partial_\eps \left(\sum_{i=1}^{m}f_i\right)(x)\subset \bigcap_{\eta>0}\bigcup_{\substack{\eps_i\ge 0
						,\;i=1,\ldots,m
						\\
						\sum_{i=1}^m\eps_i=\eps+\eta}}\sum_{i=1}^m \partial_{\eps_i}f_i(x).
			\end{gather}
			Take $l\in \partial_\eps\left(\sum_{i=1}^m f_i\right)(x)$ and $\eta>0$.
			By (iii) and Proposition~\ref{AbstractP3.2}(i), we have \begin{equation*}
				(f_1^*\square \ldots \square f^*_m)(l)=\left(\sum_{i=1}^m f_i\right)^*\le l(x)-\left(\sum_{i=1}^m f_i\right)(x)+\eps.
			\end{equation*}
			In particular, $(f_1^*\square \ldots \square f^*_m)(l)<+\infty$.
			By Definition~\ref{AbstractD3.1}(i) of the infimal convolution, there exist $l_1,\ldots,l_m\in \LL$ such that $l=\sum_{i=1}^m l_i$ and
			$$
			\sum_{i=1}^m f_i^*(l_m)< l(x)-\sum_{i=1}^m f_i(x)+\eps+\eta,
			$$
			or equivalently, $\sum_{i=1}^m\ga_i<\eps+\eta$, where
			$\gamma_i:=f_i(x)+f_i^*(l_i)-l_i(x)$ ($i=1,\ldots,m$).
			Observe that, for all $i=1,\ldots,m$, $\gamma_i\ge 0$ by Definition~\ref{AbstractD3.1}(ii) of the conjugate function. Moreover, $l_i\in \partial_{\gamma_i}f_i(x)$ by Proposition \ref{AbstractP3.2}(i).
			Set $\eps_i:= \left[ \gamma_i+(\eps+\eta-\sum_{j=1}^m\gamma_j)/m\right]$ ($i=1,\ldots,m$).
			Thus, $\eps_i>\gamma_i$, $l_i\in \partial_{\eps_i} f_i(x)$ for all $i=1,\ldots,m$, $\sum_{i=1}^m \eps_i=\eps +\eta$, and consequently,
			$$
			l\in \bigcap_{\eta>0}\,\,\bigcup_{\substack{\eps_i\ge 0
					,\;i=1,\ldots,m
					\\
					\sum_{i=1}^m\eps_i=\eps+\eta}}\sum_{i=1}^m \partial_{\eps_i}f_i(x),
			$$
			which proves \eqref{T13P2}.
			
			\item[] (iv) $\Rightarrow$ (i)
			Let $x\in X$, $\eps \ge 0$ and $K>1$.
			If $x\notin \bigcap_{i=1}^m\dom f_i(x)$, then
			$\partial_\eps \left(\sum_{i=1}^{m}f_i\right)(x)=\emptyset$, and \eqref{AbstractT3.3.1} holds automatically.
			{Assume now that $x\in \bigcap_{i=1}^m\dom f_i(x)$. Fix $\alpha>0$. We will prove that (i) holds for $K=1+\alpha$. Using (iv), we obtain
				\sloppy
				\[
				\begin{array}{rcl}
					\partial_\eps \left(\sum_{i=1}^{m}f_i\right)(x)&\subset& \bigcap_{\eta>0}\bigcup_{\substack{\eps_i\ge 0
							,\;i=1,\ldots,m
							\\
							\sum_{i=1}^m\eps_i=\eps+\eta}}\sum_{i=1}^m \partial_{\eps_i}f_i(x)\\
					&&\\
					&\subset &\bigcap_{\eta>0}\sum_{i=1}^m \partial_{\eps+\eta}f_i(x)\subset \sum_{i=1}^m
					\partial_{(1+\alpha)\eps}f_i(x),
				\end{array}
				\]
			}
			
			where we used the fact that $\eps_i\le \eps+\eta$ in the second inclusion, and $\eta:=\alpha \eps$ in the last one. The above expression { agrees with the inclusion in (i) for any $K:=1+\alpha>1$.}
		\end{enumerate}
	\end{proof}
	
	As mentioned above, Theorem~\ref{AbstractT3.3} is an extension
	to the framework of abstract convexity of the main result 
	of \cite{Bor}
	playing a key role in deriving
	constraint qualifications for zero duality gap
	in
	convex optimization.
	We quote
	here this
	result for comparison.
	
	\begin{theorem}\label{AbstractT2}
		{\rm \cite[Theorem 3.2]{Bor}}
		Let $X$ be a normed vector space, $X^*$ its conjugate space with the weak* topology,
		and $f_i: X\to \R_{+\infty}$
		$(m\ge2)$
		be proper lower semicontinuous convex functions.
		The
		following
		properties
		are equivalent:
		\begin{enumerate}
			\item
			there exists a $K>0$ such that for all $x\in \bigcap_{i=1}^m \dom f_i$ and $\eps >0$,
			\begin{equation}\label{AbstractT2.1}
				\cl{\left(\sum_{i=1}^m\partial_{\eps}f_i(x)\right)}\subset \sum_{i=1}^m \partial_{K\eps} f_i(x);
			\end{equation}
			\item
			equality \eqref{AbstractT3.3.2} hods true
			for all $x\in X$ and $\eps \ge 0$;
			\item
			$\left(\sum_{i=1}^m f_i\right)^*=f_1^*\square \ldots\square f_m^*$;
			\item
			$f_1^*\square\ldots\square f_m^*$ is 
			lower semicontinuous.
		\end{enumerate}
	\end{theorem}
	%
	%
	\begin{remark}
		{Note that inclusion \eqref{AbstractT2.1} has the same right hand side as \eqref{AbstractT3.3.1}, but
			its
			left hand side
			involves
			a topological closure expression.
			Therefore in general \eqref{AbstractT2.1} is more restrictive than \eqref{AbstractT3.3.1}.}
		We  show in
		Proposition~\ref{AbstractP3.4}
		that \eqref{AbstractT3.3.1} $\Leftrightarrow$ \eqref{AbstractT2.1}
		when
		$\LL$ is equipped with a topology possessing
		property \eqref{AbstractF2} stated below.
	\end{remark}
	
	Recall the following
	subdifferential calculus
	result for the conventional
	convex setting (cf. \cite[Corollary 2.6.7]{Zal02}),
	which is the key ingredient in the proof of the implication (i)~$\Rightarrow$~(iii) in Theorem~\ref{AbstractT2} (cf. \cite{Bor}).
	
	\begin{lemma}\label{L16}
		Let $X$ be a normed vector space, $X^*$ its conjugate space with weak* topology,
		and $f_i: X\to \R_{+\infty}$
		$(m\ge2)$
		be proper lower semicontinuous convex functions.
		Then, for all $x\in X$ and $\eps \ge 0$, 
		\begin{equation}
			\label{AbstractF2}
			\partial_{\eps}\left(\sum_{i=1}^mf_i\right)(x)
			=\bigcap_{\eta>0}
			\cl{\bigcup_{\substack{\eps_i\ge 0
						,\;i=1,\ldots,m
						\\
						\sum_{i=1}^m\eps_i=\eps+\eta}}\sum_{i=1}^m\partial_{\eps_i} f_i(x)}.
		\end{equation}
	\end{lemma}
	
	\begin{proposition}\label{AbstractP3.4}
		Let $f_1,\ldots,f_m,: X\to \R_{+\infty}$ $(m\ge 2)$, $x\in X$ and $\eps>0$.
		Suppose that $\LL$ is equipped with a topology such that
		equality \eqref{AbstractF2} holds.
		Then
		conditions \eqref{AbstractT3.3.1} and \eqref{AbstractT2.1} are equivalent.
	\end{proposition}
	
	\begin{proof}
		\if{
			Observe that for any positive numbers $\eta,\eps$ with $0<\eta \le \eps$, we have
			\begin{equation}\label{AbstractP3.4.P1}
				\bigcup_{\substack{\eps_1+\ldots+\eps_m=\eps+\eta,\\ \eps_i\ge 0}}\sum_{i=1}^m\partial_{\eps_i}f_i(x)
				\subset \sum_{i=1}^m\partial_{\eps+\eta}f_i(x)
				\subset \sum_{i=1}^m\partial_{2\eps}f_i(x).
			\end{equation}
		}\fi
		Assume that condition \eqref{AbstractT2.1} holds with some $K>0$.
		Then \eqref{AbstractT3.3.1} holds with any $K'>K$.
		Indeed,
		using \eqref{AbstractF2},
		we obtain:
		\begin{gather*}
			\partial_{\eps}\left(\sum_{i=1}^mf_i\right)(x)
			\overset{\eqref{AbstractF2}}{\subset}
			\bigcap_{\eta>0}
			\cl\left(\sum_{i=1}^m \partial_{\eps+\eta} f_i(x)\right)
			\overset{\eqref{AbstractT2.1}}{\subset}
			\bigcap_{\eta>0}
			\sum_{i=1}^m \partial_{K(\eps+\eta)}f_i(x)\subset\sum_{i=1}^m \partial_{K'\eps}f_i(x),
		\end{gather*}
		{where we used $\eta:=\alpha \eps$ for $\alpha>0$ in the last inclusion.}
		Conversely, assume that condition \eqref{AbstractT3.3.1} holds with some $K>0$.
		Then { we will show that }\eqref{AbstractT2.1} holds with $K':=mK$.
		Indeed, we start by writing 
		\[
		\bigcup_{\substack{\eps_i\ge 0
				,\;i=1,\ldots,m
				\\
				\sum_{i=1}^m\eps_i=m\eps+\eta}}\sum_{i=1}^m\partial_{\eps_i} f_i(x)\supset \sum_{i=1}^m\partial_{\eps+(\eta/m)} f_i(x)\supset \sum_{i=1}^m\partial_{\eps} f_i(x),
		\]
		where the first inclusion is obtained by using the instance $\eps_i:=\eps+(\eta/m)$ for all $i$, and the second inclusion follows from the fact that $\partial_{\eps+(\eta/m)} f_i(x)\supset \partial_{\eps} f_i(x)$.
		Taking closure on both sides of the expression above and then taking {the intersection over} all $\eta>0$ yields

		\[
		\partial_{m\eps}\left(\sum_{i=1}^mf_i\right)(x) = \bigcap_{\eta>0}
		\cl\left(\bigcup_{\substack{\eps_i\ge 0
				,\;i=1,\ldots,m
				\\
				\sum_{i=1}^m\eps_i=m\eps+\eta}}\sum_{i=1}^m\partial_{\eps_i} f_i(x)\right)\supset \cl\left(\sum_{i=1}^m\partial_{\eps}f_i(x)\right)
		\]
		
		where we used \eqref{AbstractF2} in the equality (with $m\eps$ in place of $\eps$). Altogether, we obtain
		
		\begin{gather*}
			\cl\left(\sum_{i=1}^m\partial_{\eps}f_i(x)\right)
			\overset{\eqref{AbstractF2}}{\subset} \partial_{m\eps}\left(\sum_{i=1}^mf_i\right)(x) \overset{\eqref{AbstractT3.3.1}}{\subset} \sum_{i=1}^m \partial_{Km\eps}f_i(x),
		\end{gather*}
		which gives \eqref{AbstractT2.1} for ${K'}:={mK}$ in place of $K$. This completes the proof.
	\end{proof}
	
	\begin{remark} \label{AbstractAb_con_R1}
		\begin{enumerate}
			\item
			From the proof above we see that, if  \eqref{AbstractT2.1} holds with $K>0$, then \eqref{AbstractT3.3.1} holds with
			any $K'>K$.
			On the other hand, if \eqref{AbstractT3.3.1} holds with $K>0$, then  \eqref{AbstractT2.1} holds with $mK$.
			\item
			When $0\in \LL$, property (iii) in Theorem~\ref{AbstractT3.3}
			(property (iii)
			in Theorem~\ref{AbstractT2})
			ensures
			condition \eqref{AbstractEx-Zero}, and consequently,
			the zero duality gap
			property.
			In
			Theorem~\ref{AbstractT3.3},
			the functions $f_1,\ldots,f_m$ are
			not assumed to be
			convex.
			However, without convexity, condition \eqref{AbstractT3.3.1} may not be easy to satisfy. When
			for some $i=1,\ldots,m$,
			the function $f_i$
			is
			not convex,
			in view of Proposition~\ref{AbstractP4.1}(iii),
			there exists a point $x\in X$
			such that $\partial_\eps f_i(x)=\emptyset$
			for some $\eps >0$.
			In this situation, condition \eqref{AbstractT3.3.1} fails to hold.
		\end{enumerate}
	\end{remark}
	
	The zero duality gap property is equivalent to equality \eqref{AbstractEx-Zero}, which is  clearly less restrictive than property (iii) in Theorem~\ref{AbstractT3.3}.
	In the next theorem, we relax condition \eqref{AbstractT3.3.1} as well as property (iii) in Theorem~\ref{AbstractT3.3},
	and
	obtain a characterization of the zero duality gap property, which is new even in the classical convex case.
	
	
	
	
	\begin{theorem}\label{AbstractT3.5}
		Let $f_1,\ldots,f_m,: X\to \R_{+\infty}$ $(m\ge 2)$ be such that $\bigcap_{i=1}^m\dom f_i\neq \emptyset$.  Suppose that $0\in \LL$.
		The following properties are equivalent:
		\begin{enumerate}
			\item
			$0\in\bigcap_{\eps>0}\left(\sum_{i=1}^m\partial_{\eps}f_i\right)(X)$;
			
			\item  
			$\left(\sum_{i=1}^m f_i\right)^*(0)=f_1^*\square\ldots\square f_m^*(0)<+\infty$.
		\end{enumerate}
	\end{theorem}
	
	\begin{proof}
		\begin{enumerate}
			\item[] (i) $\Rightarrow$ (ii).
			In view of Proposition~\ref{AbstractP3.2}(iv), we only need to show that
			\begin{equation}\label{T19P1}
				f_1^*\square\ldots\square f_m^*(0)\le
				\left(\sum_{i=1}^m f_i\right)^*(0)<+\infty.
			\end{equation}
			Let $\eps >0$.
			By (i), there exists an $x\in X$ such that $0\in\sum_{i=1}^m \partial_{\eps/m}f_i(x)$.
			Hence, in view of Proposition~\ref{AbstractP3.2}(v),
			$0\in\partial_{\eps}\left(\sum_{i=1}^mf_i\right)(x)$.
			Since $\eps >0$ is arbitrary, it follows from Proposition~\ref{AbstractP3.2}(ii) that $0\in  \dom \left(\sum_{i=1}^m f_i\right)^*$, or equivalently, $\left(\sum_{i=1}^m f_i\right)^*(0)<+\infty$.
			{ Recall that $0\in\sum_{i=1}^m \partial_{\eps/m}f_i(x)$, hence }
			there exist $l_i\in \partial_{\eps/m}f_i(x)$ ($i=1,\ldots,m$) such that $\sum_{i=1}^ml_i=0$ and
			$
			f_i^*(l_i)+f_i(x)\le l_i(x)+\eps/m$ for all $i=1,\ldots,m.
			$
			Hence, the definitions of the infimal convolution and conjugate function yield
			$$
			f_1^*\square\ldots\square f_m^*(0)\le \sum_{i=1}^m f_i^*(l_i)
			\le \eps-\sum_{i=1}^mf_i(x)
			{\le}\left(\sum_{i=1}^m f_i\right)^*(0)+\eps,
			$$
			Letting $\eps\downarrow 0$, we arrive at the first inequality in \eqref{T19P1}.
			
			\item[] (ii) $\Rightarrow$ (i).
			Let $\eps >0$.
			{By (ii), $0\in \dom \left(\sum_{i=1}^m f_i\right)^*$. The latter fact, together with Proposition~\ref{AbstractP3.2}(ii), imply that}
			there exists $x\in X$ such that
			$
			0\in \partial_{2\eps}\left(\sum_{i=1}^m f_i\right)(x).
			$
			Hence, by (ii) and Proposition~\ref{AbstractP3.2}(i),
			\begin{equation}
				\label{AbstractT3.5.P1}
				f_1^*\square\ldots\square f_m^*(0)
				=\left(\sum_{i=1}^m f_i\right)^*(0)
				{\le}2\eps-\left(\sum_{i=1}^m f_i\right)(x).
			\end{equation}
			By Definition~\ref{AbstractD3.1}(i) of the infimal convolution, there are $l_1,\ldots,l_m\in\LL$ such that \begin{equation}\label{AbstractT5.5.P2}
				\sum_{i=1}^ml_i=0\quad\text{and}\quad
				\sum_{i=1}^mf_i^*(l_i)
				<f_1^*\square\ldots\square f_m^*(0)+\eps.
			\end{equation}
			Combining \eqref{AbstractT3.5.P1} and \eqref{AbstractT5.5.P2}, we obtain
			\begin{gather}\label{AbstractT5.5.P3}
				\sum_{i=1}^m(f_i^*(l_i)+f_i(x)-l_i(x))<\eps.
			\end{gather}
			{By Definition~\ref{AbstractD3.1}(ii) of the conjugate function, we know that $f_i(x)+f_i^*(l_i)-l_i(x)\ge 0$ for all $i=1,\ldots,m$. Hence, each term in \eqref{AbstractT5.5.P3} must be less than $\eps$. Namely, } for all $i=1,\ldots,m$,
			$f_i^*(l_i)+f_i(x)<l_i(x)+\eps$,
			and consequently, by Proposition~\ref{AbstractP3.2}(i),
			$l_i\in \partial_{\eps}f_i(x)$. We conclude that
			$0=\sum_{i=1}^ml_i\in  \sum_{i=1}^m \partial_{\eps}f_i(x){\subset \bigcup_{x'\in X} \sum_{i=1}^m \partial_{\eps}f_i(x')= (\sum_{i=1}^m\partial_{\eps}f_i)(X)}$. {Since $\eps>0$ is arbitrary, (i) is established.}
			
		\end{enumerate}
	\end{proof}
	
	The next theorem provides another pair of equivalent properties, each of which is stronger than the corresponding property in Theorem~\ref{AbstractT3.5}.
	
	\begin{theorem}\label{AbstractT3.6}
		Let $f_1,\ldots,f_m,: X\to \R_{+\infty}$ $(m\ge 2)$ and $x\in\bigcap_{i=1}^m\dom f_i$.  Suppose that $0\in \LL$.
		The following properties are equivalent:
		\begin{enumerate}
			\item
			$0\in\bigcap_{\eps>0}\sum_{i=1}^m\partial_{\eps}f_i(x)$;
			\item  
			$\left(\sum_{i=1}^m f_i\right)^*(0)=f_1^*\square\ldots\square f_m^*(0)<+\infty$, and $x$ is a solution of problem \eqref{AbstractP}.
		\end{enumerate}
	\end{theorem}
	
	\begin{proof}
		\begin{enumerate}
			\item[] (i) $\Rightarrow$ (ii).
			The first assertion in (ii) follows from (i) thanks to Theorem~\ref{AbstractT3.5}.
			Fix $\eps>0$.
			By (i), there are $l_1,\ldots,l_m\in \LL$ such that $\sum_{i=1}^m l_i=0$, and $l_{i}\in \partial_{\eps} f_i(x)$ for all $i=1,\ldots,m$, and consequently, by Proposition~\ref{AbstractP3.2}(i), $f_i^*(l_i)+f_i(x)\le  l_i(x)+\eps$ for all $i=1,\ldots,m$.
			Hence, recalling the definitions of the infimal convolution and conjugate function, we obtain
			$$
			f_1^*\square\ldots\square f_m^*(0)
			{\le} \sum_{i=1}^m f_i^*(l_i)\le -\left(\sum_{i=1}^m f_i\right)(x)+m\eps{\le} \left(\sum_{i=1}^m f_i\right)^*(0)+m\eps.
			$$
			Letting $\eps\downarrow 0$, and taking into account Proposition~\ref{AbstractP3.2}(iv), we get
			\begin{equation}\label{AbstractT3.6.P1}
				f_1^*\square\ldots\square f_m^*(0)=-\left(\sum_{i=1}^m f_i\right)(x)= \left(\sum_{i=1}^m f_i\right)^*(0).
			\end{equation}
			By \eqref{AbstractP1}, we deduce that $v(P)=\sum_{i=1}^m f_i(x)$.
			
			\item[] (ii) $\Rightarrow$ (i).
			By the assumptions in (ii), we have \eqref{AbstractT3.6.P1}.
			Fix $\eps>0$. By Definition~\ref{AbstractD3.1}(i) of the infimal convolution, there are $l_{1},\ldots,l_{m}\in \LL$
			satisfying conditions \eqref{AbstractT5.5.P2}.
			Combining \eqref{AbstractT5.5.P2} and \eqref{AbstractT3.6.P1}, we arrive at \eqref{AbstractT5.5.P3}.
			Now we proceed as in the proof of Theorem~\ref{AbstractT3.5},
			{ to conclude that $0=\sum_{i=1}^ml_i\in  \sum_{i=1}^m \partial_{\eps}f_i(x)$. Since $\eps>0$ is arbitrary, (i) is established.}
		\end{enumerate}
	\end{proof}
	
	
	\section{Abstract Convexity with Pointwise Convergence Topology}\label{S4}
	
	In this section, we
	assume the space $\LL$ of abstract linear functions to be a vector space equipped with
	a certain
	topology
	and
	expand some fundamental results of standard convex analysis to the framework of abstract convexity. Condition \eqref{AbstractT2.1},
	playing a key role in deriving constraint qualifications for
	zero duality gap in convex optimization in \cite{Bor},
	is fully extended into the abstract convexity framework. 
	
	\subsection{Pointwise Convergence Topology}
	
	
	In accordance with Definition \ref{AbstractALS}, the space $\LL$ of abstract linear functions possesses the standard addition operation.
	We are now going to assume additionally that $\LL$ is closed with respect to the natural scalar multiplication: given an $l\in\LL$ and $\al\in\R$, the function $\al l$ defined by
	$\al l(x):=\alpha l(x)$ for all $x\in X$, belongs to $\LL$.
	This makes $\LL$ a vector space.
	As a consequence, $\HH$ is also a vector space, and they both are subspaces of the ambient vector space $\mathcal{F}$ of all real-valued functions on $X$.
	
	From now on, we assume that $\LL$, $\HH$ and $\mathcal{F}$ are all equipped with the pointwise convergence topology, i.e.
	the weakest topology that makes all functions $f\mapsto f(x)$,
	$x\in X$,
	continuous (on respective spaces).
	Such a topology is often referred to as the weak* topology.
	Recall that
	no topology is assumed on $X$.
	If $X$ is a normed vector space and $\LL$ is the space of all linear functions on $X$, then the pointwise convergence topology on $\LL$ discussed above coincides with the conventional weak* topology.
	
	The two assertions in the next proposition are natural generalizations of the corresponding facts from
	\cite{Bre}.
	Their proofs require only cosmetic modifications.
	
	\begin{proposition}\label{AbstractP4.11}
		\begin{enumerate}
			\item
			Let
			$l_0\in \LL$, $\eps >0$, $n\in\N$,
			and $x_1,\ldots,x_n\in X$.
			Then the set
			$$
			V(x_1,\ldots,x_n|\eps)=\{l\in \LL: |l(x_i)-l_0(x_i)|<\eps,\quad \forall i=1,\ldots,n \}
			$$
			is a neighbourhood of $l_0$ in $\LL$.
			Moreover, the collection of all $V(x_1,\ldots,x_n|\eps)$, $n\in \N$, $x_1,\ldots,x_n\in X$ and $\eps >0$ forms a basis of neighbourhoods of $l_0$ in $\LL$.
			\item
			If
			$\LL$ is 
			closed in $\mathcal{F}$, then it is a locally convex topological vector space.
		\end{enumerate}
	\end{proposition}
	
	\begin{example}
		\cite[Example~2.1]{Kay}
		\label{AbstractE4.3}
		Let $\LL$ be a set of functions defined on the Euclidean space $\R^n$, comprising all
		functions $l:=\sum_{i=0}^{n}a^ih_i\in \LL$ where $a^i\in \R$, $i=0,1,\ldots,n$, and
		$$
		h_0(x)=\norm{x}^2,\quad h_1(x)=x_1,\ldots,h_n(x)=x_n \quad\text{for all}\quad x=(x_1,\ldots,x_n)\in \R^n.
		$$
		Obviously $\LL$ is a vector space.
		Moreover, considered with
		the pointwise convergence topology,
		it possesses the following properties.
		\begin{enumerate}
			\item
			$\LL$ is 
			closed in $\mathcal{F}$.
			\item
			Given a net $(l_i)_{i\in I}\subset\LL$ and an element $l\in\LL$ with
			$l_i:=a_i^0h_0+\ldots+a_i^nh_n$ and $l:=a^0h_0+\ldots+a^nh_n$, one has $l_i\to l$
			if and only if $a_i^j\to a^j$
			for all $j=0,1,\ldots,n$.
			\item
			$\LL$
			is {\it homeomorphic} to $\R^{n+1}$ with the standard Euclidean norm.
		\end{enumerate}
		
		\begin{proof}
			\begin{enumerate}
				\item
				Let a net $(l_i)_{i\in I}\subset\LL$ with $l_i:=a_i^0 h_0+\ldots+a_i^nh_n$ converge to some $l\in \mathcal F$.
				We are going to show that $l\in\LL$.
				Take the unit vector $e_0:=(1,0,\ldots,0)\in\R^{n}$.
				We have
				$l_i(e_0)+l_i(-e_0)=2a^0_i\to l(e_0)+l(-e_0)=:2a^0$.
				Thus,
				$a_i^0\to a^0$,
				and consequently,
				$l_i(x)-a_i^0\norm{x}^2\to l(x)-a^0\norm{x}^2$ for all $x\in X$.
				{For every $i\in I$,
					$l_i-a_i^0\norm{\cdot}^2=a_i^1 h_1+\ldots+a_i^nh_n$, and, hence, it is a usual linear function on $\R^n$.}
				Thus, $l-a^0\norm{\cdot}^2$ is also a linear function on $\R^n$
				as
				a pointwise limit of
				linear functions.
				This implies that $l-a^0\norm{\cdot}^2
				=\sum_{i=1}^{n}a^ih_i$ for some $a^i\in \R$, $i=1,\ldots,n$.
				Therefore,  $l\in \LL$.
				
				\item
				If $l_i\to l$, then, as shown in (i), we have that $a_i^0\to a^0$ and $l_i-a_i^0\norm{\cdot}^2\to l-a^0\norm{\cdot}^2$. {Taking $x:=e^j$ the canonical vectors in $\R^n$ for  all $j=1,\ldots,n$, we deduce from the latter limit that $l_i(e^j)-a_i^0=a_i^j-a_i^0\to a^j-a^0$. Using also (i), we obtain that $a^j_i\to a^j$ for all $j=0,\ldots,n$. The converse statement is trivial.}
				
				\item
				Consider the bijective mapping $\varphi: \LL\to \R^{n+1}$,
				given for any
				$l=\sum_{i=0}^n a^ih_i\in\LL$ by
				$\varphi(l):=(a^0,a^1,\ldots,a^n)$.
				Observe that $\varphi^{-1}(a^0,a^1,\ldots,a^n)=\sum_{j=0}^na^jh_j$.
				In view of (ii), both $\varphi$ and $\varphi^{-1}$ are continuous.
			\end{enumerate}
		\end{proof}
	\end{example}
	
	\if{
		In view of Example~\ref{AbstractE4.3}, we can treat the space of abstract linear functions $\LL$ (in Example~\ref{AbstractE4.3}) as the vector space $\R^{n+1}$ with the basis $h_0, h_1,\ldots,h_n$. Moreover, Example~\ref{AbstractE4.3}(iii) implies that the weak$^*$ topology is equivalent to the Euclidean norm topology in $\R^{n+1}$. 
	}\fi
	
	
	The weak$^*$ compactness of the unit ball $ \B^*$ in the space dual to a normed vector space
	is of utmost importance
	in analysis and applications.
	Here we establish a generalization of the Banach--Alaoglu--Bourbaki theorem (see \cite[Theorem 3.16]{Bre}) to
	the space $\mathcal{F}$ of real-valued functions on $X$.
	
	
	
	\begin{theorem}\label{AbstractT4.3}
		Let $F,G\in\FF$ and
		$G\le F$.
		Then the set $A:=\{f\in \mathcal{F}: G\le f\le F\}$ is 
		compact.
		%
	\end{theorem}
	
	\begin{proof}
		Given a function $f\in\FF$, set $\phi(f):=(f(x))_{x\in X}$.
		Then $\phi$ is a homeomorphism between $\mathcal{F}$ and $\R^X$ considered with the product topology.
		By Tychonoff's theorem, the set
		$H=\prod_{x\in X}[G(x),F(x)]$
		is compact in $\R^X$ as a product of compact sets.
		Hence, $A=\phi\iv(H)$ is 
		compact.
	\end{proof}
	
	
	
	
	\begin{corollary}\label{AbstractC4.4}
		If a subset $A\subset\mathcal{F}$ is 
		closed, and the functions $F:=\sup_{f\in A}f$ and $G:=\inf_{f\in A}f$ are everywhere finite, then $A$ is
		compact.
	\end{corollary}
	
	\begin{proof}
		Under the assumptions, $F,G\in\FF$,
		$G\le F$ and $A\subset\{f\in \mathcal{F}: G\le f\le F\}$, and consequently, $A=\{f\in \mathcal{F}: G\le f\le F\}\cap A$.
		As $A$ is 
		closed, the conclusion follows from Theorem~\ref{AbstractT4.3}.
	\end{proof}
	
	\begin{remark}
		If $X$ is a normed vector space, $\LL=X^*$ and $A=\{x^*: \norm{x^*}\le 1\}$, Corollary~\ref{AbstractC4.4} { with $F(x):=\sup_{x^*\in A} \langle x^*,x \rangle=\|x\|$ and  $G(x):=\inf_{x^*\in A} \langle x^*,x \rangle=-\|x\|$} recaptures the Banach--Alaoglu--Bourbaki theorem.
	\end{remark}
	
	
	%

	\subsection{Sum Rule for Subdifferentials}
	
	In this subsection, we
	establish a version of
	the following (extended) sum rule for abstract convex functions
	under a weaker qualification condition.
	
	\begin{theorem}
		\label{AbstractKay} {\rm \cite[Theorem 3.2]{Kay}}
		Let $f,g: X\to \R_{+\infty}$ be $\HH-$convex functions with $\dom f \cap \dom g \neq \emptyset$.
		The following conditions are equivalent:
		\begin{gather}
			\label{Abstractcon*}
			\epi f^*+ \epi g^*= \epi(f+g)^*,
			\\
			\label{Abstractkay1}
			\partial_{\eps}(f+g)(x)=\bigcup_{\eps_1,\eps_2 \ge0,\; \eps_1+\eps_2=\eps}\big(\partial_{\eps_1} f(x)+\partial_{\eps_2} g(x)\big)
		\end{gather}
		for all $x\in \dom f \cap \dom g$ and $\eps \ge 0$.
	\end{theorem}
	
	The additivity property \eqref{Abstractcon*} of the epigraphs of the conjugate functions is used in Theorem~\ref{AbstractKay} to ensure the sum rule \eqref{Abstractkay1}. This assumption is rather strong, since it entails the 
	closedness of the set in the \LHS\ of \eqref{Abstractcon*}, which is not true in general even for conventional convex functions.
	Condition \eqref{Abstractcon*} is stronger than the widely used (e.g. in Theorem~\ref{AbstractT2}) qualification condition $(f+g)^*=f^*\square g^*$; cf. \cite[Corollary 5.1]{Kay}. We are going to replace condition \eqref{Abstractcon*} by the following less restrictive additivity property:
	\begin{equation}\label{Abstractcon}
		\cl(\epi f^*+\epi g^*)= \epi(f+g)^*,
	\end{equation}
	where $\cl$ stands for the 
	closure in $\LL\times \R$.
	The 
	topology on $\LL\times \R$ is the product of the
	pointwise convergence
	topology on $\LL$ 
	adopted here
	and the standard topology on $\R$.
	Condition \eqref{Abstractcon} always holds for conventional convex functions; cf. \cite[Corollary 3.1]{Kay} and \cite{Str19.2}.
	
	We next extend some classical facts of convex analysis to the framework of abstract convexity.
	
	\begin{proposition}\label{AbstractP2.P1}
		\begin{enumerate}
			\item
			For any $x\in X$, the function $\LL\ni l\mapsto\phi_x(l):=l(x)$ is 
			continuous;
			\item
			for any function $f:X\to \R_{+\infty}$, its conjugate function $f^*$ is 
			lower semicontinuous;
			\item
			$\LL-$convex sets are
			closed;
			\if{
				\AK{22/03/20.
					What about $\HH-$convex sets?}
				
				\HB{If we replace the set $\LL$ by $\HH$ and if $\HH$ is weak* closed, then (iii) is applicable for $\HH-$convex sets.}
			}\fi
			
			\item
			$\HH-$convex sets are
			closed;
			
			\item
			for any $f: X \to \R_{+\infty}$, $x\in \dom f$ and $\eps \ge 0$,
			the $\eps$-subdifferential $\partial_\eps f(x)$ is 
			closed;
			
			\item
			for any functions $f_1,\ldots,f_m$ $(m\ge 2)$, $x\in \bigcap_{i=1}^m\dom f_i$, and $\eps \ge 0$,
			\begin{equation*}
				\cl
				{\bigcup_{\substack{\eps_i\ge 0
							,\;i=1,\ldots,n
							\\
							\sum_{i=1}^m\eps_i=\eps}} \sum_{i=1}^m\partial_{\eps_i}f_i(x)}
				\subset
				\partial_{\eps}\left(\sum_{i=1}^mf_i\right)(x),
			\end{equation*}
			and consequently,
			\begin{equation}\label{AbstractP4.6.1}
				\bigcap_{\eta>0}\cl
				{\bigcup_{\substack{\eps_i\ge 0
							,\;i=1,\ldots,n
							\\
							\sum_{i=1}^m\eps_i=\eps+\eta}} \sum_{i=1}^m\partial_{\eps_i}f_i(x)}
				\subset
				\partial_{\eps}\left(\sum_{i=1}^mf_i\right)(x).
			\end{equation}
		\end{enumerate}
	\end{proposition}
	
	\begin{proof}
		\begin{enumerate}
			\item
			{Recall that we are considering the pointwise convergence in $\LL$. Let $x\in X$ and take any net $(l_i)_{i\in I}\subset\LL$ converging to some $l\in\LL$.} Then
			$
			\phi_x(l_i)=l_i(x)\to l(x)=\phi_x(l),
			$
			i.e. $\phi_x$ is 
			continuous.
			\item
			{We will show that $\epi f^*$ is closed in $\LL\times \R$.} Indeed,
			Let $f:X\to \R_{+\infty}$. Then,
			$$
			\epi f^*=\{(l,\la): f^*(l)\le \la\}=\bigcap_{x\in X}\{(l,\la): l(x)-f(x))\le \la \}=\bigcap_{x\in X} \epi (\phi_x-f(x)).
			$$
			{Using (i) we have, that for each fixed $x$,} the function $\phi_x-f(x): \LL \to \R_{+\infty}$ is 
			continuous, and consequently, the set $\epi(\phi_x-f(x))$ is 
			closed in $\LL\times \R$. Hence, so is $\epi f^*$. This implies the 
			lower semicontinuity of the function $f^*$.
			
			\item
			Let $\es\ne A\in\LL$ be $\LL-$convex.
			By  Proposition~\ref{AbstractP4.2.2}(i), there is an $\LL-$convex function $f:X\to\R_{+\infty}$ such that
			$
			A=\supp f.
			$
			Consider the set $B:=\cl A$
			and the function $g:=\sup_{l\in B}l$.
			Obviously $f\le g$.
			Take any $x\in X$ and $l\in B$.
			Then there exists a net $(l_i)_{i\in I}\subset A$ converging to $l$ in the 
			pointwise convergence
			topology.
			We have $l_i(x)\le f(x)$ for all $i\in I$, and consequently, $l(x)\le f(x)$.
			Since $l\in B$ is arbitrary, it follows that $g(x)\le f(x)$, and since $x\in X$ is arbitrary, it follows that $g\le f$.
			Consequently, $g=f$.
			Thus, if $l\in B$, then $l\le f$, i.e. $l\in A$, and consequently, $B=A$.
			
			\item
			Thanks to Remark~\ref{R7}(iii), the assertion is a consequence of (iii).
			
			\item
			Let $f: X \to \R_{+\infty}$, $x\in \dom f$ and $\eps \ge 0$.
			By Proposition~\ref{AbstractP3.2}(i),
			$
			\partial_{\eps} f(x)=\{l\in\LL: f^*(l)+f(x)\le l(x)+\eps\}=\{l\in\LL: f^*(l)-l(x)\le \eps-f(x)\}.
			$ Namely, it is a level set of the function $l\mapsto f^*(l)-l(x)$.
			By (i) and (ii), the latter function is 
			lower semicontinuous. Indeed,  $f^*$ is 
			lower semicontinuous by (ii) and $l\mapsto l(x)$ is
			continuous by (i).
			Hence, $\partial_{\eps} f(x)$ is 
			closed.
			
			\item
			Thanks to (v), the right-hand side in the first inclusion is closed. This fact, together with Proposition~\ref{AbstractP3.2}(v), readily gives the first inclusion.
			The second inclusion is a consequence of the first one written with $\eps+\eta$ in place of $\eps$, {followed by an application of } Proposition~\ref{AbstractP3.2}(iii).
		\end{enumerate}
	\end{proof}

	The next theorem presents an extension of an important result of classical convex analysis; cf. \cite{LiNg08}.
	
	\begin{theorem}\label{AbstractSum}
		Let
		$f_1,\ldots,f_m:X\to \R_{+\infty}$ $(m\ge 2)$
		and
		$\cap_{i=1}^m \dom f_i\neq \emptyset$. Assume that
		\begin{equation}\label{AbstractCon}
			\cl\left(\sum_{i=1}^m\epi f_i^*\right)=\epi \left(\sum_{i=1}^m f_i\right)^* .
		\end{equation}
		Then the sum rule \eqref{AbstractF2} holds for all $x\in \cap_{i=1}^m \dom f_i$ and $\eps \ge 0$.
	\end{theorem}
	
	\begin{proof}
		Let $x\in \cap_{i=1}^m \dom f_i$ and $\eps \ge 0$.
		By Proposition~\ref{AbstractP2.P1}(vi), inclusion \eqref{AbstractP4.6.1} holds.
		Now we prove the converse inclusion.
		Let $l\in \partial_{\eps}\left(\sum_{j=1}^mf_j\right)(x)$.
		By Proposition~\ref{AbstractP3.2}(i), $\left(\sum_{j=1}^mf_j\right)^*(l) +\left(\sum_{j=1}^mf_j\right)(x)\le l(x)+\eps$.
		Thus, making use of \eqref{AbstractCon}, we have
		$$\left(l,l(x)+\eps-\sum_{j=1}^mf_j(x)\right)\in\epi \left(\sum_{j=1}^mf_j\right)^*{=}\cl\left(\sum_{j=1}^m \epi f_j^* \right).$$
		There are nets $(l_{1,i}, \la_{1,i})_{i\in I}, \ldots, (l_{m,i}, \la_{m,i})_{i\in I}$
		such that
		\begin{gather}
			\label{AbstractAb-sum.P1}
			f_j^*(l_{j,i})\le\la_{j,i},\quad j=1,\ldots,m,\; i\in I,\\
			\label{AbstractAb-sum.P2}
			\sum_{j=1}^m l_{j,i}{\to} l,\quad  \sum_{j=1}^m \la_{j,i}{\to}l(x)+\eps-\sum_{j=1}^mf_j(x).
		\end{gather}
		For all $j=1,\ldots,m$, $i\in I$, set $\gamma_{j,i}:=\la_{j,i}+f_j(x)-l_{j,i}(x)$.
		Thanks to Proposition~\ref{AbstractP4.1}(iv), inequalities \eqref{AbstractAb-sum.P1} imply $\gamma_{j,i}\ge 0$ and $l_{j,i}\in \partial_{\gamma_{j,i}}f_j(x)$.
		By \eqref{AbstractAb-sum.P2}, $\sum_{j=1}^m\gamma_{j,i}\to\eps$. Since each $\gamma_{j,i}\ge 0$, we can assume, without loss of generality, $\gamma_{j,i}\to\ga_j'\ge0$, $j=1,\ldots,m$, and $\sum_{j=1}^m\ga_j'=\eps$.
		Fix $\eta >0$.
		There is an $i_0\in I$ such that, for all $i\ge i_0$ and $j=1,\ldots,m$, we have
		$\gamma_{j,i}<\eps_j:=\ga_j' +\eta/m$.
		Thus, $l_{j,i}\in \partial_{\eps_{j}}f_j(x)$, $j=1,\ldots,m$, and $\sum_{j=1}^m\eps_{j}=\eps+\eta$.
		It follows that
		$$l\in
		\cl{\bigcup_{\substack{\eps_i\ge 0
					,\;i=1,\ldots,n
					\\
					\sum_{i=1}^m\eps_i=\eps+\eta}}\sum_{i=1}^m\partial_{\eps_i} f_i(x)}.$$
		Since $\eta>0$ is arbitrary, the proof is complete.
	\end{proof}
	
	\begin{remark}
		The exact sum rule \eqref{Abstractkay1} in Theorem~\ref{AbstractKay} is a stronger property than the sum rule \eqref{AbstractF2} in Lemma~\ref{L16} and Theorem~\ref{AbstractSum}.
	\end{remark}
	
	\subsection{Zero Duality Gap with Pointwise Convergence Topology}
	
	We next present the main theorem of this section.
	
	\begin{theorem}\label{AbstractT3}
		Let
		$f_1,\ldots,f_m:X\to \R_{+\infty}$ $(m\ge 2)$.
		Consider the following
		properties:
		\begin{enumerate}
			\item
			there exists a $K>1$ such that condition \eqref{AbstractT3.3.1} holds
			for all $x\in \bigcap_{i=1}^m \dom f_i$ and  $\eps>0$;
			\item
			there exists a $K>0$ such that condition \eqref{AbstractT2.1} holds
			for all $x\in \bigcap_{i=1}^m \dom f_i$ and  $\eps>0$;
			\item
			condition
			\eqref{AbstractT3.3.2} holds
			for all $x\in X$ and $\eps \ge 0$;
			\item
			$\left(\sum_{i=1}^m f_i\right)^*=f_1^*\square \ldots\square f_m^*$;
			\item
			$f_1^*\square\ldots\square f_m^*$ is 
			lower semicontinuous.
		\end{enumerate}
		We have {\rm (i) $\Leftrightarrow$ (iii) $\Leftrightarrow$ (iv)  $\Rightarrow$ (v) $\Rightarrow$ (ii)}.
		If the sum rule \eqref{AbstractF2} (or condition \eqref{AbstractCon}) holds, then all five
		properties
		are equivalent.
	\end{theorem}
	
	\begin{proof}
		The equivalence of (i), (iii) and (iv) is shown in Theorem~\ref{AbstractT3.3}.
		Property (v) is a consequence of (iv) since $\left(\sum_{i=1}^n f_i\right)^*$ is 
		lower semicontinuous by Proposition~\ref{AbstractP2.P1}(ii).
		{Even though not trivial, the proof of the implication (v) $\Rightarrow$ (ii) proceeds exactly the same way (mutatis mutandis) as in the proof of (iv) $\Rightarrow$ (i) in Theorem~\ref{AbstractT2}; cf.
			\cite[Theorem 3.2]{Bor}. Hence, we omit the proof here.}
		\if{
			\AK{21/03/20.
				To be checked.
				I have removed the assumption that $X$ is a normed vector space.
				Could it be needed here?}
			\HB{I think so. I will check.}
		}\fi
		If \eqref{AbstractF2} holds, then by  Proposition~\ref{AbstractP3.4}, (i) $\Leftrightarrow$ (ii) { and all the statements become equivalent.}
	\end{proof}
	
	\begin{remark}
		Similarly to \cite{Bor},
		Theorem~\ref{AbstractT3} entails a series of important corollaries with proofs almost identical to those of the corresponding statements in \cite{Bor}.
	\end{remark}
	
	\section{Zero Duality Gap for a Family of Nonconvex Problems}
	\label{AbstractExample_AC}
	
	In this section, we consider a nontrivial
	abstract convexity framework
	in which the 
	additivity condition \eqref{Abstractcon},
	and as a consequence also the sum rule \eqref{AbstractF2}
	hold for all $\HH-$convex functions.
	Throughout this section,
	$\FF$ is the family of all functions from $\R$ to $\R$.
	The spaces $\LL$ of abstract linear and $\HH$ of abstract affine functions are defined as follows.
	
	
	\begin{definition}\label{AbstractEX.D1}
		Given $a,b\in\R$, set $\phi_{a}(x):=ax^2$ and $\psi_{a,b}(x):=ax^2+b$ for all $x\in \R$.
		Next set
		$\LL:=\{\phi_a: a\in \R\}$ and $\HH:=\{\psi_{a,b}:a,b\in \R\}$.
	\end{definition}
	
	Clearly, $\psi_{a,b}=\phi_{a}+b$ and $\phi_{a}=\psi_{a,0}$.
	The sets $\LL$ and $\HH$ defined above obviously satisfy all the conditions in Definitions \ref{AbstractALS} and \ref{AbstractD4.1}, respectively.
	Moreover, they are vector spaces.
	
	\begin{remark}\label{AbstractR3}
		Similarly to
		Example~\ref{AbstractE4.3},
		it is easy to show
		that the 
		pointwise convergence
		topologies on $\LL$ and $\HH$ are isomorphic to the usual topologies of $\R$ and $\R^2$, respectively.
		The mappings $\LL\ni\phi_{a}\mapsto a$ and $\HH\ni\psi_{a,b}\mapsto(a,b)$ are homeomorphisms. 
	\end{remark}
	
	{First, we characterize $\LL-$convex functions and $\LL-$convex sets.}
	
	\begin{proposition}[$\LL-$convex functions and sets]
		\label{AbstractP4.E0}
		\begin{enumerate}
			\item
			A function $f:\R\to \R_{+\infty}$ is
			$\LL-$convex if and only if
			either $f\in\LL$ or
			$$
			f(x)=
			\begin{cases}
				0& \text{if } x=0,\\
				+\infty& \text{otherwise}.
			\end{cases}
			$$
			\item
			A nonempty subset $C\subset \LL$ is $\LL-$convex if and only if either $C=\{\phi_a: a\le\bar a\}$ for some $\bar a\in\R$ or $C=\LL$.
		\end{enumerate}
	\end{proposition}
	
	\begin{proof}
		\begin{enumerate}
			\item
			By Definition~\ref{AbstractD4.2}(ii),
			a function $f:\R\to \R_{+\infty}$ is
			$\LL-$convex if and only if $f=\sup_{a\in A}\phi_a$ for some subset $A\subset\R$.
			Since $f>-\infty$, we deduce that $A\ne\es$. Set $\bar a:=\sup A$.
			If $\bar a<+\infty$, then $f=\phi_{\bar a}\in\LL$. Indeed, since $a\le \bar a$ for all $a\in A$, the definition of $f$ yields $f\le \phi_{\bar a}$.
			For the opposite inequality, take any $\eta>0$.  There exists $a_{\eta}\in A$ such that $\bar a -\eta <a_{\eta}\le \bar a$. Using the definition of $f$ again, we have that $f(x)\ge  a_{\eta}x^2 > (\bar a -\eta)x^2$ for every $x$. Since this is true for  
			every $\eta>0$, we deduce that $f\ge \phi_{\bar a}$. Altogether, we have shown that  $f=\phi_{\bar a}\in \LL$ if $\sup A=\bar a<+\infty$. Otherwise, $A$ is unbounded above and, hence, by definition $f(x)=+\infty$ whenever $x\ne0$ and $f(0)=0$.
			\item
			By Proposition~\ref{AbstractP4.2.2}(i), a nonempty subset $C\subset \LL$ is $\LL-$convex if and only if $C=\supp f$ for some $\LL-$convex function $f:\R\to \R_{+\infty}$.
			In view of (i), either there exists an $\bar a\in\R$ such that $f=\phi_{\bar a}$, and consequently,
			$C=\supp\phi_{\bar a}=\{\phi_a: \phi_a\le\phi_{\bar a}\}=\{\phi_a: a\le\bar a\}$, or $C=\{\phi_a: \phi_a(x)<+\infty \text{ for all } x\ne0\}=\{\phi_a: a\in\R\}=\LL$.
		\end{enumerate}
	\end{proof}
	
	{Denote by $\R^2_+:=\{(a,b)\in\R^2:a\ge0,\;b\ge0\}$. For $A,B\subset \R^2$, define \\ \[A-B:=\{x-y\::\: x\in A,\, y\in B\}.\]
		Next, we study $\HH-$convex functions and $\HH-$convex sets.}

	\begin{proposition}[$\HH-$convex functions and sets]
		\label{AbstractP4.E1}
		\begin{enumerate}
			\item
			An $\HH-$convex function is even and \lsc.
			
			\item
			Assume that $C$ is $\HH-$convex, and let $A\subset \R^2$ be such that $C=\{\psi_{a,b}\in \HH:(a,b)\in A\}$. Then $A$ is closed and convex, and $A-\R^2_+\subset A$.
			\item
			Assume that $A\subset \R^2$ is
			closed and convex, and such that it also verifies that $A-\R^2_+\subset A$ and 
			$\sup_{(a,b)\in A}\psi_{a,b}\not\equiv+\infty$. Then the induced set 
			$C:=\{\psi_{a,b}\in \HH:(a,b)\in A\}$ is $\HH-$convex. 
		\end{enumerate}
	\end{proposition}
	
	\begin{proof}
		\begin{enumerate}
			\item
			By Definition~\ref{AbstractD4.2}(iii),
			if a function $f:\R\to \R_{+\infty}$ is
			$\HH-$convex, then there exists a subset $A\subset\R^2$ such that
			$f=\sup_{(a,b)\in A}\psi_{a,b}$.
			The assertion follows since all functions $\psi_{a,b}$ are even and continuous.
			\item
			Take $A$ and $C$ as in (ii).  By Proposition~\ref{AbstractP4.2.2}(iii), $C$ is convex in the conventional sense. Let $(a_1,b_1),(a_2,b_2)\in A$, $\al\in[0,1]$ and $(\hat a,\hat b) :=\al(a_1,b_1)+(1-\al)(a_2,b_2)$.
			Then, in view of Definition~\ref{AbstractEX.D1}, $\psi_{\hat a,\hat b} =\al\psi_{a_1,b_1}+(1-\al)\psi_{a_2,b_2}\in C$, i.e. $(\hat a,\hat b)\in A$, and consequently, $A$ is convex.
			By Proposition~\ref{AbstractP2.P1}(iv), $C$ is 
			closed.
			Let a sequence $(a_k,b_k)\in A$ converge to some $(\bar a,\bar b)\in\R^2$.
			Then, for any $x\in\R$, $\psi_{a_k,b_k}(x)\to\psi_{\bar a,\bar b}(x)$.
			{Since $C$ is closed,} $\psi_{\bar a,\bar b}\in C$, i.e. $(\bar a,\bar b)\in A$, and consequently, $A$ is closed.
			By Proposition~\ref{AbstractP4.2.2}(ii),
			$C=\supp f$ for some
			$\HH-$convex function $f:\R\to\R_{+\infty}$.
			Let $(a,b)\in A$, and
			$a'\le a$, $b'\le b$.
			Then
			$\psi_{a',b'}\le \psi_{a,b}\le f$.
			Thus, $\psi_{a',b'}\in \supp f=C$, and consequently, $(a',b')\in A$.
			Hence, $A-\R^2_+\subset A$.
			\if{
				Let $\bar b=\infty$.
				Then $\bar a:=\sup_{(a,b)\in A}a<\infty$ since otherwise $C=\HH$.
				Choose a number $a_0>\bar a$.
				We have $\psi_{a_0,0}\notin C$ and, for any $t\in\R$, $\psi_{a_0,0}(t)=a_0 t^2<\sup_{\psi\in C}\psi(t)=+\infty$, i.e. $C$ is not $\HH-$convex.
				A contradiction.
				Hence, $\bar b<\infty$.
			}\fi
			
			\item
			Suppose that $A$ satisfies the conditions in the assertion, and
			take $\psi_{a_0,b_0}\notin C$ for some $a_0,b_0\in\R$. { By definition of $A$, this means that $(a_0,b_0)\notin A$. Define $\bar a:=\sup_{(a,b)\in A}a\le +\infty$. We will consider two cases: either $a_0>\bar a$ or $a_0\le \bar a$. Assume that $a_0>\bar a$.} Since $\sup_{(a,b)\in A} \psi_{a,b}\not\equiv+\infty$, there exist numbers $x_0\ge0$ and $\mu\in \R$ such that $\psi_{a,b}(x_0)=ax_0^2+b<\mu$ for all $(a,b)\in A$.
			For a sufficiently large $x>x_0$, we have $(a_0-\bar a) x^2>\mu-\bar ax_0^2-b_0$.
			Hence, for any $(a,b)\in A$,
			\begin{equation*}
				\psi_{a_0,b_0}(x)=a_0x^2+b_0>\bar a(x^2-x_0^2)+\mu\ge a(x^2-x_0^2)+\mu>ax^2+b=\psi_{a,b}(x),
			\end{equation*}
			and consequently, $\psi_{a_0,b_0}(x)>\sup_{\psi\in C}\psi(x)$. By definition, this implies that $C$ is $\HH-$convex.
			
			{ Assume now that $a_0\le\bar a$. Since $A$ is closed and convex, and $(a_0,b_0)\notin A$, we can apply} the strict convex separation theorem (Lemma~\ref{AbstractS.G.F}), to obtain a nonzero vector $(\al,\be)\in \R^2$ such that
			\begin{equation}\label{P36P2}
				\al a_0+\be b_0> \sup_{(a,b)\in A}(\al a+\be b).
			\end{equation}
			Since $A-\R^2_+\subset A$, we have $\al\ge0$ and $\be\ge0$.
			Moreover, $\be\ne0$ because $a_0\le\bar a$.
			Thus, $\be>0$.
			Set $x:=\sqrt{\al/\be}$.
			Inequality \eqref{P36P2} yields
			$$\psi_{a_0,b_0}(x)=a_0x^2+b_0 >\sup_{(a,b)\in A}(ax^2+b)=\sup_{\psi\in C}\psi(x).$$
			By Definition~\ref{AbstractD4.2}(v), this implies that $C$ is $\HH-$con\-vex. The proof is complete.
		\end{enumerate}
	\end{proof}
	
	\begin{remark}
		The properties in Proposition~\ref{AbstractP4.E1}(i) are necessary but not sufficient. Indeed, the function $x\mapsto\cos x$ is even and continuous but not $\HH-$convex.
	\end{remark}
	The assumption that $\sup_{(a,b)\in A}\psi_{a,b}\not\equiv+\infty$ in Proposition~\ref{AbstractP4.E1}(iii) is essential, as the next example shows.
	
	\begin{example}
		Let $C:=\{\psi_{a,b}:a\le0,\,b\in\R\}$.
		We have $\psi_{1,0}\notin C$ and, for any $x\in\R$, $\psi_{1,0}(x)=x^2<\sup_{\psi\in C}\psi(x)=+\infty$. By Definition~\ref{AbstractD4.2}(v), this implies that $C$ is not $\HH-$convex.
		Observe that
		$C$ is a proper subset of $\HH$.
		The whole space $\HH$ is obviously $\HH-$convex; thus, the assumption $\sup_{(a,b)\in A}\psi_{a,b}\not\equiv+\infty$ {is in general not necessary for the set to be $\HH-$convex.}
	\end{example}
	
	In the above example $\sup_{(a,b)\in A}b =+\infty$. The next example shows that this condition does not necessarily entail the absence of $\HH-$convexity.
	
	\begin{example}
		Let $C:=\{\psi_{a,b}:a+b\le0\}$.
		We have $\sup_{(a,b)\in A}a=\sup_{(a,b)\in A}b =+\infty$.
		At the same time, the set $A:=\{(a,b):a+b\le0\}$ satisfies the conditions in Proposition~\ref{AbstractP4.E1}(iii).
		In particular, $\sup_{(a,b)\in A}\psi_{a,b}(1)=0<+\infty$.
		Hence, $C$ is $\HH-$convex.
	\end{example}
	
	\if{
		\AK{27/04/20.
			Currently there is a gap between the necessary and sufficient conditions.
			This leads to problems in Proposition~\ref{AbstractAb-Ex}.
			Could the sufficient conditions be weakened having the above example in mind?}
	}\fi
	
	\if{
		\AK{29/03/20.
			I am afraid obtaining descriptions of $\HH-$convex functions and sets is not that simple.
			$\HH$ is a two-parameter family, and we are facing a vector optimization problem.
			Considering the parameters separately as in  \eqref{AbstractDf.A} and \eqref{AbstractDf.B} is unlikely to lead to a meaningful description.
			Besides, once the $+\infty$ value is allowed, we cannot talk about `upper boundedness' as in part (i) below.
			
			`Abstract' convexity and weak$^*$ closedness should also be reformulated in terms of $a$ and $b$ as in Definition~\ref{AbstractEX.D1}.}
	}\fi
	
	In accordance with Proposition~\ref{AbstractP4.E1}(iii), any set of the form $\{\psi_{a,b}: a\le\bar a,b\le\bar b\}$ with some $\bar a\in \R\cup\{+\infty\}$ and $\bar b\in \R$ is $\HH-$convex.
	However, not every $\HH-$convex set is of this form.
	
\begin{figure}[!ht]
		\begin{minipage}[t]{.5\linewidth}
			\centering
			\includegraphics{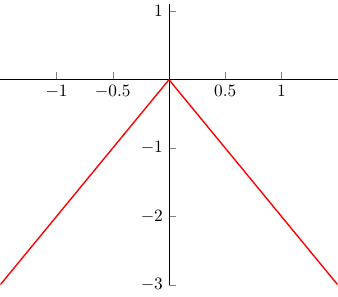}
		\end{minipage}
		\begin{minipage}[t]{.5\linewidth}
			\centering
			\includegraphics{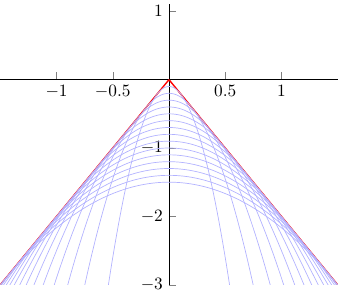}
		\end{minipage}
		\caption{Function $-2\abs{x}$ and its
			support set}
		\label{FF1}
	\end{figure}

	\begin{example}
		The set $C:=\{\psi_{a,b}: a<0 ,\; ab\ge1\}$ (the support set of the function $x\mapsto-2{\abs{x}}$; cf. Fig.~\ref{FF1}) satisfies the conditions in Proposition~\ref{AbstractP4.E1}(iii), and consequently is $\HH-$convex, but is not of the form $\{\psi_{a,b}: a\le\bar a,b\le\bar b\}$.
	\end{example}
	
	The next statement shows that $\HH-$convex functions in Definition~\ref{AbstractEX.D1} satisfy the additivity property \eqref{Abstractcon}.
	
	\begin{proposition}\label{AbstractAb-Ex}
		Suppose $f,g:\R\to\R_{+\infty}$ are $\HH-$convex and {$\dom f \cap \dom g \neq \emptyset$}. Then
		$\cl(\supp f+\supp g)$ is $\HH-$convex. {Moreover, we have that}
		\sloppy
		\begin{equation}\label{T40-1}
			\supp(f+g)=\cl(\supp f+\supp g).
		\end{equation}
		{Equivalently,} condition \eqref{Abstractcon} is satisfied.
	\end{proposition}
	
	\begin{proof}
		{By assumption, $f+g\not\equiv+\infty$, and hence $f\not\equiv+\infty$ and $g\not\equiv+\infty$.
			Using 
			Remark \ref{Abstract-R}(i) and Proposition~\ref{AbstractP4.2.2}(ii), the sets
			$\supp f$ and $\supp g$ are non-empty and $\HH-$convex.}
		{Let $A_f,A_g\subset\R^2$ be such that $\supp f=\{\psi_{a,b}:(a,b)\in A_f\}$ and $\supp g=\{\psi_{a,b}:(a,b)\in A_g\}$. By Proposition~\ref{AbstractP4.E1}(ii),
			$A_f$ and $A_g$ are closed, convex, and such that $A_f-\R^2_+\subset A_f$, $A_g-\R^2_+\subset A_g$.
			To prove that the set $A:=\cl(A_f+A_g)$ is $\HH-$convex, we will show that it satisfies the properties in Proposition~\ref{AbstractP4.E1}(ii).} Indeed, it is obviously closed and convex. Moreover, since $2\R^2_+=\R^2_+$, we have
		\begin{equation}\label{A2}
			(A_f+A_g)-\R^2_+=(A_f-\R^2_+)+(A_g-\R^2_+)\subset A_f+A_g\subset A.
		\end{equation}
		Since $A$ is closed, we have that $A-\R^2_+\subset A$. 
		{Indeed, let us first show that
			$$A - \R^2=\cl(A_f +A_g)-\R^2\subset \cl((A_f+A_g)-\R^2),$$
			where the first equality follows from the definition of $A$. To prove the inclusion, take $x\in \cl(A_f+A_g)-\R^2$, or $x = a - b$ with $a\in \cl(A_f+A_g)$ and $b \in \R^2$. There are $a_n \in A_f+A_g$ ($n\ge 0$) such that $a_n \to a$. Since $a_n - b \in (A_f+A_g)-\R^2$ and $a_n - b\to (a - b) = x$, we have $x\in \cl((A_f+A_g)-\R^2)$. Then,
			$$
			A - \R^2=\cl(A_f+A_g)-\R^2\subset \cl((A_f+A_g)-\R^2)\subset A,
			$$
			where the last inclusion follows by taking closures in \eqref{A2} and using the fact that $A$ is closed. Hence, our claim holds and $A-\R^2_+\subset A$.} {Since $\dom f \cap \dom g \neq \emptyset$, we have $\sup_{(a,b)\in A}\psi_{a,b} = f+g \not\equiv +\infty$. By Proposition~\ref{AbstractP4.E1}(iii), this implies $A$ is a $\HH-$convex set. The first assertion in the proposition is established.} {We will show next that $f+g$ is $\HH-$convex. Indeed, taking into account the closedness of $A$ again, $$\sup_{(a,b)\in A}\psi_{a,b}=\sup_{(a,b)\in A_f+A_g}\psi_{a,b}=\sup_{(a_f,b_f)\in A_f} \psi_{a_f,b_f}+\sup_{(a_g,b_g)\in A_g}\psi_{a_g,b_g}=f+g,$$}
		{where the second equality follows from the separability of the variables, and the fact that for every  $(a_f,b_f)\in A_f$ and $(a_g,b_g)\in A_g$, we have $\psi_{a,b}=\psi_{a_f,b_f}+\psi_{a_g,b_g}$, where $a=a_f+a_g$ and $b=b_f+b_g$. The above equality shows that $f+g$ is $\HH-$convex by Definition \ref{AbstractD4.2}(iii). Since $A$ is $\HH-$convex. Remark \ref{Abstract-R}(iv) implies that 
			\[
			\{\psi_{a,b}\::\: (a,b)\in A\}=\supp(f+g).
			\]
			From the definitions, {one can check  directly} that 
			\[
			(a,b)\in A \hbox{ if and only if }\psi_{a,b}\in \cl(\supp\, f +\supp\, g).
			\]
			The two expressions above yield  \eqref{T40-1}.} Recall from Definition~\ref{AbstractEX.D1} that $\psi_{a,b}=\phi_a+b$ for all $a,b\in\R$. In view of Proposition~\ref{AbstractP4.1}(ii), $\psi_{a,b}\in\supp f$ if and only if $(\phi_a,-b)\in\epi f^*$, and similar equivalences hold for the support sets of $g$ and $f+g$ and the epigraphs of their conjugate functions.
		Hence, condition \eqref{T40-1} can be rewritten equivalently as \eqref{Abstractcon}.
	\end{proof}
	
\begin{figure}[!ht]
		\begin{minipage}[t]{.5\linewidth}
			\centering
			\includegraphics[width=0.8\textwidth]{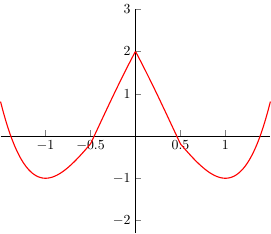}
			\caption{$f_1+f_2+f_3$.}
			\label{abstract_Fig1}
		\end{minipage}
		\begin{minipage}[t]{.5\linewidth}
			\centering
			\includegraphics[width=1.0\textwidth]{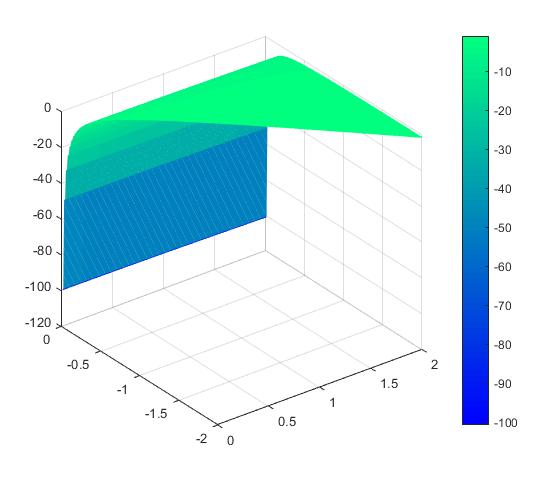}
			\caption{Restriction of $-f_1^*-f_2^*-f_3^*$ to the subspace ${\phi_{a_1}+\phi_{a_2}+\phi_{a_3}=0}$.}
			\label{abstract_Fig2}
		\end{minipage}
	\end{figure}

	Now as an illustration we consider an example of the minimization problem \eqref{AbstractP} and its corresponding dual problem \eqref{AbstractD}:
	\begin{gather}
		\label{p}\tag{p}
		\inf\left(f_1(x)+f_2(x)+f_3(x)\right)\quad
		\textnormal{s.t.}\quad x\in\R,
		\\
		\label{d}\tag{d}
		\sup\left(-f_1^*(\phi_{a_1})-f_2^*(\phi_{a_2})-f_3^*(\phi_{a_3})\right) \quad
		\textnormal{s.t.}\quad \phi_{a_1},\phi_{a_2},\phi_{a_3}\in \LL,\;\;
		\phi_{a_1}+\phi_{a_2}+\phi_{a_3}=0,
	\end{gather}
	where, for all $x\in\R$,
	\begin{align}
		\label{AbstractP2-f_1}
		f_1(x):=&x^4-x^2,\quad
		f_2(x):=1-2\abs{x},\\
		\label{AbstractP2-f_3}
		f_3(x):=&
		\begin{cases}
			1-2\abs{x}& \text{if }-\frac{1}{2}\le x\le\frac{1}{2},\\
			0& \text{otherwise}.
		\end{cases}
	\end{align}
	The objective functions of problems \eqref{p} and \eqref{d} are presented in Fig.~\ref{abstract_Fig1} and Fig.~\ref{abstract_Fig2}, respectively.
	It can be seen from Fig.~\ref{abstract_Fig1} that $x=\pm 1$ are the global solutions of problem \eqref{p}, with optimal value equal $-1$.
	The objective function of problem \eqref{d} is concave.
	(This is always the case because a conjugate function is the supremum of functions which are linear in $\LL$.)
	Hence, problem \eqref{d} is a conventional convex programming problem (of maximizing a concave function over a subspace).
	
	The next three propositions provide representations of the conjugates, support sets and subdifferentials of the functions $f_1, f_2, f_3$.
	The second one also states that the functions are $\HH-$convex.
	Their proofs are given in the Appendix.
	
	\begin{proposition}\label{Conj}
		The conjugates of
		the functions $f_1, f_2, f_3$ given
		by \eqref{AbstractP2-f_1} and \eqref{AbstractP2-f_3},
		have the following representations:
		\begin{align}
			\label{AbstractP2-f_1*}
			f^*_1(\phi_a)=&
			\begin{cases}
				\frac{(a+1)^2}{4}& \text{if } a\ge -1,\\
				0& \text{if }a<-1,
			\end{cases}\\
			\label{AbstractP2-f_2*}
			f_2^*(\phi_a)=&
			\begin{cases}
				+\infty& \text{if } a\ge0,\\
				-1-\frac{1}{a}& \text{if }a<0,
			\end{cases}\\
			\label{AbstractP2-f_3*}
			f_3^*(\phi_a)=&
			\begin{cases}
				+\infty& \text{if } a>0,\\
				\frac{a}{4}& \text{if } -2\le a\le0,\\
				-1-\frac{1}{a}& \text{if }a<-2.
			\end{cases}
		\end{align}
	\end{proposition}
	
\begin{figure}[!ht]
		\begin{minipage}[t]{.5\linewidth}
			\centering
			\includegraphics[width=0.8\textwidth]{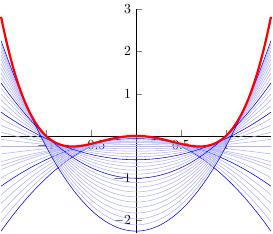}
			\caption{Support set of $f_1$.}
			\label{Fig4}
		\end{minipage}
		\begin{minipage}[t]{.5\linewidth}
			\centering
			\includegraphics[width=0.8\textwidth]{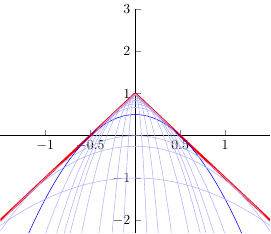}
			\caption{Support set of $f_2$.}
			\label{Fig5}
		\end{minipage}
	\end{figure}

\begin{figure}[!ht]
		\begin{minipage}[t]{.5\linewidth}
			\centering
			\includegraphics[width=0.8\textwidth]{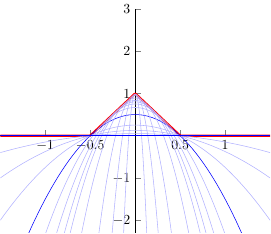}
			\caption{Support set of $f_3$.}
			\label{Fig6}
		\end{minipage}
		\begin{minipage}[t]{.5\linewidth}
			\centering
			\includegraphics[width=0.8\textwidth]{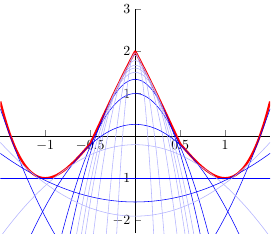}
			\caption{Support set of $f_1+f_2+f_3$.}
		\end{minipage}
	\end{figure}

	\begin{proposition}\label{AbstractP44}
		The functions $f_1, f_2, f_3$ given
		by \eqref{AbstractP2-f_1} and \eqref{AbstractP2-f_3}
		are $\HH-$convex.
		Their support sets  have the following representations (cf. Fig.~\ref{Fig4}, \ref{Fig5} and \ref{Fig6}):
		\begin{align}
			\label{AbstractE4.1.}
			\supp f_1=&\left\{\psi_{a,b}: \frac{(a+1)^2}{4}+b\le0\right\}\bigcup \left\{\psi_{a,b}:a\le -1,\; b\le 0\right\},\\
			\label{AbstractE4.2.1}
			\supp f_2=&\left\{\psi_{a,b}: a<0,\; b\le\frac{1}{a}+1\right\},\\
			\label{AbstractE4.3.1}
			\supp f_3=&\left\{\psi_{a,b}: a\le-2,\; b\le\frac{1}{a}+1\right\}\bigcup\left\{\psi_{a,b}: -2<a\le0,\;b\le-\frac{a}{4}\right\}.
		\end{align}
	\end{proposition}
	
	\begin{proposition}\label{Subdifferential}
		The subdifferentials of
		the functions $f_1, f_2, f_3$ given
		by \eqref{AbstractP2-f_1} and \eqref{AbstractP2-f_3},
		have the following representations:
		\begin{align}
			\label{AbstractE4.1.2}
			\partial f_1(x)=&\begin{cases}
				\{\phi_{a}: a\le -1\}&\text{if } x=0,\\
				\{\phi_{a}: a=2x^2-1\}& \text{otherwise},
			\end{cases}
			\\
			\label{AbstractE4.2.2}
			\partial f_2(x)=&\begin{cases}
				\emptyset&\text{if } x =0,\\
				\left\{\phi_{-1/{\abs{x}}}\right\}& \text{otherwise},
			\end{cases}
			\\\label{AbstractE4.3.2}
			\partial f_3(x)=&\begin{cases}
				\emptyset&\text{if } x =0,\\
				\left\{\phi_{-1/{\abs{x}}}\right\} &\text{if } 0<|x|<\frac{1}{2},\\
				\left\{\phi_{a}: -2\le a\le0\right\}&\text{if } x=\pm\frac{1}{2},\\
				\{0\}& \text{otherwise}.
			\end{cases}
		\end{align}
	\end{proposition}
	
	We next show that 
	zero duality gap holds for problems \eqref{p} and \eqref{d}.
	
	\begin{proposition}
		$\left(f_1+f_2+f_3\right)^*(0)=f_1^*\square f_2^* \square f_3^*(0)<+\infty$, and consequently, zero duality gap holds for problems \eqref{p} and \eqref{d}.
		Moreover, the points $x=\pm 1$ are solutions of problem~\eqref{p}.
	\end{proposition}
	
	\begin{proof}
		From Proposition~\ref{Subdifferential},
		$\partial f_1 (\pm1)=\{\phi_1\}$,
		$\partial f_2 (\pm1)=\{\phi_{-1}\}$ and
		$\partial f_3 (\pm1)=\{0\}$.
		Hence,
		$\partial f_1 (\pm1)+ \partial f_2(\pm1)+\partial f_3(\pm1)=\{0\}$.
		The conclusion follows from Theorem~\ref{AbstractT3.6}.
	\end{proof}
	
	
	
	\section*{Appendix: Proofs of Propositions~\ref{Conj}, \ref{AbstractP44} and \ref{Subdifferential}}
	\label{apendix}
	
	\subsubsection*{Proof of Proposition~\ref{Conj}}
	
	\paragraph*{Expression~\eqref{AbstractP2-f_1*}.}\quad
	Given $a\in\R$, we have
	$f^*_1(\phi_a)=\sup_{x\in \R}\left\{(a+1)x^2-x^4\right\}$.
	If $a<-1$, then $(a+1)x^2-x^4\le0$, and the maximum (equal $0$) is attained at $x=0$; hence $f^*_1(\phi_a)=0$.
	If $a\ge-1$, then $(a+1)x^2-x^4 =\frac{(a+1)^2}{4}-\left(x^2-\frac{a+1}{2}\right)^2$, and the maximum of the expression is attained when $x^2=\frac{a+1}{2}$; hence $f^*_1(\phi_a)=\frac{(a+1)^2}{4}$.
	
	\paragraph*{Expression~\eqref{AbstractP2-f_2*}.}\quad
	Given $a\in\R$, we have
	$f^*_2(\phi_a)=\sup_{x\in \R}\left\{ax^2+2|x|\right\}-1$.
	If $a\ge0$, then obviously $f^*_2(\phi_a)=+\infty$.
	If $a<0$, then $ax^2+2|x|={a}(|x|+\frac{1}{a})^2-\frac{1}{a}$, and the maximum of the expression is attained when $|x|=-\frac{1}{a}$; hence $f^*_2(\phi_a)=-\frac{1}{a}-1$.
	
	\paragraph*{Expression~\eqref{AbstractP2-f_3*}.}\quad
	Given $a\in\R$, we have
	$f^*_3(\phi_a)=\sup_{x\in \R}\left\{ax^2+\min\{2|x|-1,0\}\right\}$.
	If $a>0$, then obviously $f^*_3(\phi_a)=+\infty$.
	Otherwise,
	$
	f^*_3(\phi_a)=\sup_{|x|\le\frac{1}{2}}\left\{ax^2+2|x|\right\}-1.
	$
	If $a<-2$, then, as shown above,
	$f^*_3(\phi_a)=-\frac{1}{a}-1$, while if $-2\le a\le0$, $f^*_3(\phi_a)=\frac{a}{4}+2\cdot\frac{1}{2}-1=\frac{a}{4}$.

	\subsubsection*{Proof of Proposition~\ref{AbstractP44}}
	\paragraph*{Expression~\eqref{AbstractE4.1.} and  $\HH-$convexity of $f_1$.}\quad
	If $\frac{(a+1)^2}{4}+b\le0$,
	then, for all $x\in \R$,
	$$
	f_1(x)=x^4-x^2=\left(x^2-\frac{a+1}{2}\right)^2-\left(\frac{(a+1)^2}{4}+b\right) +(ax^2+b)\ge ax^2+b=\psi_{a,b}(x),
	$$
	i.e. $\psi_{a,b}\in \supp f_1$.
	Similarly, if $a\le -1$ and $b\le 0$,
	then, for all $x\in \R$,
	$$
	f_1(x)=x^4-x^2=x^4-\left((a+1)x^2+b\right)+(ax^2+b)\ge ax^2+b=\psi_{a,b}(x),
	$$
	i.e. $\left\{\psi_{a,b}:a\le -1,\; b\le 0\right\}\subset \supp f_1$.
	Thus,
	\begin{align}\label{P41P1}
		\left\{\psi_{a,b}: \frac{(a+1)^2}{4}+b\le0\right\}\bigcup \left\{\psi_{a,b}:a\le -1,\; b\le 0\right\}\subset \supp f_1.
	\end{align}
	{Conversely, for any $\psi_{a,b} \in \supp f_1$, we have two possibilities: either $a\le -1$ or $a>-1$. In the first case, we claim that we must have $b\le 0$. Indeed, if $b>0$, then $f_1(0)=0<b=\psi_{a,b}(0)$; hence, $\psi_{a,b} \notin \supp f_1$, contradiction. Hence,
		$\supp f_1\subset \{\psi_{a,b}\::\: b\le 0,\, a\le -1\}$. If $a>-1$ then, with $x:=\sqrt{\frac{a+1}{2}}$, we have
		$$f_1(x)=\frac{(a+1)^2}{4}-\frac{a+1}{2}{\ge} a(\frac{a+1}{2})+b=\psi_{a,b}(x);$$
		which re-arranges as $\frac{(a+1)^2}{4}+b\le 0$. Hence, $\supp f_1\subset\{\psi_{a,b}\::\: 
		\frac{(a+1)^2}{4}+b\le 0\}$.} It follows that the inclusion in \eqref{P41P1} holds as equality.
	This proves \eqref{AbstractE4.1.}. {Let us prove now the $\HH-$convexity of $f_1$.} Given any $x\in \R$, set $a:= 2x^2-1$ and $b:=-x^4$.
	Then, $\frac{(a+1)^2}{4}+b=0$ and, in view of \eqref{AbstractE4.1.}, $\psi_{a,b}\in \supp f_1$.
	Moreover, $\psi_{a,b}(x)= ax^2+b= (2x^2-1)x^2-x^4= x^4-x^2=f_1(x)$.
	Thus,
	$
	\sup_{\psi\in \supp f_1}=f_1.
	$
	By Proposition~\ref{AbstractP4.1}(i), $f_1$ is $\HH-$convex.
	
	\paragraph*{Expression~\eqref{AbstractE4.2.1} and  $\HH-$convexity of $f_2$.}\quad
	
	\medskip
	
	Let $a<0$ and $b\le\frac{1}{a}+1$. We claim that $\psi_{a,b}\in \supp f_2$. Indeed, it is elementary to show that $\frac{1}{a}=\min_{x\in \R}\{-2|x|-ax^2\}$. By definition of $b$, for any $x\in\R$ we have
	$b\le \frac{1}{a}+1\le 1-2\abs{x}-ax^2$, and consequently, $\psi_{a,b}(x)=ax^2+b\le 1-2\abs{x}=f_2(x)$.
	Thus,
	\begin{align}\label{P41P3}
		\left\{\psi_{a,b}: a<0,\; b\le\frac{1}{a}+1\right\}\subset\supp f_2.
	\end{align}
	If $a>0$, then $\lim_{x\to \infty} \psi_{a,b}(x)=+\infty$ while $f_2\le1$, and consequently, $\psi_{a,b}(x)> f_2(x)$ for some $x\in\R$.
	Similarly, if $b>\frac{1}{a}+1$, then there is an $x\in\R$ such that $b>1-2\abs{x}-ax^2$, and consequently, $\psi_{a,b}(x)=ax^2+b>1-2\abs{x}=f_2(x)$.
	Hence, $\psi_{a,b} \notin \supp f_2$.
	It follows that the inclusion in \eqref{P41P3} holds as equality.
	This proves \eqref{AbstractE4.2.1}.
	
	Next we show that
	$\sup_{\psi\in \supp f_2}=f_2$.
	Given any $x\in\R\setminus \{0\}$, set $a:=-\frac{1}{\abs{x}}$ and $b:=\frac{1}{a}+1=1-|x|$.
	In view of \eqref{AbstractE4.2.1}, $\psi_{a,b}\in\supp f_2$.
	Moreover,
	$\psi_{a,b}(x)=ax^2+b=1-2\abs{x}=f_2(x)$.
	For $x=0$ and any $n\ge1$, set $a_n:=-n$ and $b_n:=\frac{1}{a_n}+1=1-\frac{1}{n}$. In view of \eqref{AbstractE4.2.1}, $\psi_{a_n,b_n}\in\supp f_2$.
	Moreover, $\lim_{n\to \infty}\psi_{a_n,b_n}(0)=\lim_{n\to \infty}\left(1-\frac{1}{n}\right)=1=f_2(0)$.
	Thus,
	$\sup_{\psi\in \supp f_2}=f_2$.
	By Proposition~\ref{AbstractP4.1}(i), $f_2$ is $\HH-$convex.
	
	\paragraph*{Expression~\eqref{AbstractE4.3.1} and  $\HH-$convexity of $f_3$.}\quad
	For all $a\in \R$, define 
	\[
	\chi(a):=\max_{|x|\le\frac{1}{2}}\left\{ax^2 +2|x|\right\}-1.
	\]
	As shown in the proof of Proposition~\ref{Conj},
	$\chi(a)=
	-\frac{1}{a}-1$ if $a<-2$, and
	$\chi(a)=\frac{a}{4}$ if $-2\le a\le0$.
	Thus, the \RHS\ of \eqref{AbstractE4.3.1} can be rewritten as $\left\{\psi_{a,b}: a\le 0,\; b\le-\chi(a)\right\}$.
	
	Let $a\le 0$ and $b\le-\chi(a)$. We claim that $\psi_{a,b}\in \supp f_3$.
	Indeed, if
	$|x|\le\frac{1}{2}$, then $b\le-\chi(a)\le 1-2\abs{x}-ax^2$, and consequently, $\psi_{a,b}(x)=ax^2+b\le 1-2\abs{x}=f_3(x)$. In particular, the latter fact yields 
	\begin{equation}\label{psi1}
		\psi_{a,b}(\frac{1}{2})\le 1-2\abs{1/2}= 0.
	\end{equation}
	If $|x|>\frac{1}{2}$, since $a\le 0$, we have 
	\[
	\psi_{a,b}(x)\le \psi_{a,b}(\frac{1}{2})\le  0=f_3(\frac{1}{2})=f_3(x),
	\]
	where we used \eqref{psi1} in the second inequality and the definition of $f_3$ in the equalities. Thus,
	\begin{align}\label{P41P2}
		\left\{\psi_{a,b}: a\le 0,\; b\le -\chi_a(x)\right\}\subset\supp f_3.
	\end{align}
	If $a>0$, then $\lim_{x\to \infty} \psi_{a,b}(x)=+\infty$ while $f_3\le1$, and consequently, $\psi_{a,b}(x)> f_3(x)$ for some $x\in\R$.
	Similarly, if $b>-\chi(a)$, then there is an $x\in [-\frac{1}{2},\frac{1}{2}]$ such that $b>1-2\abs{x}-ax^2$, and consequently, $\psi_{a,b}(x)=ax^2+b>1-2\abs{x}=f_3(x)$.
	Hence, $\psi_{a,b} \notin \supp f_3$.
	It follows that the inclusion in \eqref{P41P2} holds as equality.
	This proves \eqref{AbstractE4.3.1}.
	
	Next we show that
	$\sup_{\psi\in \supp f_3}=f_3$.
	Given any $x\notin(-\frac{1}{2},\frac{1}{2})$, set $a:=0$ and $b:=-\chi(a)=0$.
	In view of \eqref{AbstractE4.3.1}, $\psi_{a,b}\in\supp f_3$.
	Hence, we have that 
	$\psi_{a,b}(x)=f_3(x)=0$.
	Given any $x\in(-\frac{1}{2},\frac{1}{2})\setminus \{0\}$, set $a:=-\frac{1}{\abs{x}}\in(-\infty,-2]$ and $b:=-\chi(a)=\frac{1}{a}+1=1-|x|$.
	In view of \eqref{AbstractE4.3.1}, $\psi_{a,b}\in\supp f_3$.
	Moreover,
	$\psi_{a,b}(x)=ax^2+b=1-2\abs{x}=f_3(x)$.
	For $x=0$ and any $n>1$, set $a_n:=-n$ and $b_n:=-\chi(a_n)=\frac{1}{a_n}+1=1-\frac{1}{n}$. In view of \eqref{AbstractE4.3.1}, $\psi_{a_n,b_n}\in\supp f_3$.
	Moreover, $\lim_{n\to \infty}\psi_{a_n,b_n}(0)=\lim_{n\to \infty}\left(1-\frac{1}{n}\right)=1=f_3(0)$.
	Thus,
	$\sup_{\psi\in \supp f_3}=f_3$.
	By Proposition~\ref{AbstractP4.1}(i), $f_3$ is $\HH-$convex.

	\subsubsection*{Proof of Proposition~\ref{Subdifferential}}
	\paragraph*{Expression~\eqref{AbstractE4.1.2}.}\quad
	By Proposition~\ref{AbstractP3.2}(i), $\phi_{a}\in\partial f_1(x)$ for some $x,a\in\R$ if and only if
	$f_1^*(\phi_{a})+f_1(x)\le\phi_{a}(x)$,
	or equivalently, 
	\begin{gather}
		\label{P42P1}
		f_1^*(\phi_{a})\le(a+1)x^2-x^4.
	\end{gather}
	If $x=0$, inequality \eqref{P42P1} becomes
	\[
	f_1^*(\phi_{a})\le 0,
	\]
	which, by \eqref{AbstractP2-f_1*} happens if and only if $a\le-1$; hence $\partial f_1(0)=\{\phi_{a}: a\le -1\}$. Let $x\ne0$. By \eqref{AbstractP2-f_1*}, if $a\le-1$, then \eqref{P42P1} reduces to $0\le(a+1)x^2-x^4$, which cannot be satisfied by a nonzero $x$.
	Thus, if  $x\ne0$, \eqref{P42P1} can only happen if $a>-1$, which implies, by \eqref{AbstractP2-f_1*}
	\[
	\frac{(a+1)^2}{4}\le (a+1)x^2- x^4,
	\]
	or $\left(x^2-\frac{a+1}{2} \right)^2\le0$, i.e. $a=2x^2-1$; hence $\partial f_1(x)=\{\phi_{a}: a=2x^2-1\}$.
	
	\paragraph*{Expression~\eqref{AbstractE4.2.2}.}\quad
	By Proposition~\ref{AbstractP3.2}(i), $\phi_{a}\in\partial f_2(x)$ for some $x,a\in\R$ if and only if
	$f_2^*(\phi_{a})+f_2(x)\le\phi_{a}(x)$,
	or equivalently,
	\begin{gather}
		\label{P42P2}
		f_2^*(\phi_{a})\le ax^2+2|x|-1.
	\end{gather}
	Recall the expression for $f_2^*$ given by \eqref{AbstractP2-f_2*}.
	Inequality \eqref{P42P2} implies $a<0$. If $x=0$,
	inequality \eqref{P42P2} yields $f_2^*(\phi_{a})\le-1$.
	However, in view of \eqref{AbstractP2-f_2*},
	$f_2^*(\phi_{a})>-1$ for all $a\in\R$.
	Hence, $\partial f_2(0)=\es$.
	Let $x\ne0$ and
	$a<0$. By \eqref{AbstractP2-f_2*}, $f_2^*(\phi_{a})=-1-\frac{1}{a}$, and \eqref{P42P2} becomes $ax^2+2|x|+\frac{1}{a}\ge0$, {or equivalently, $a(|x|+\frac{1}{a})^2\ge0$. Since $a<0$, this translates as
		$(|x|+\frac{1}{a})^2\le0$,} i.e. $a=-\frac{1}{\abs{x}}$. Hence, $\partial f_2(x) =\left\{\phi_{-1/{\abs{x}}}\right\}$.
	
	\paragraph*{Expression~\eqref{AbstractE4.3.2}.}\quad
	By Proposition~\ref{AbstractP3.2}(i), $\phi_{a}\in\partial f_3(x)$ for some $x,a\in\R$ if and only if
	$f_3^*(\phi_{a})+f_3(x)\le\phi_{a}(x)$,
	or equivalently,
	\begin{gather}
		\label{P42P3}
		f_3^*(\phi_{a})\le ax^2-f_3(x).
	\end{gather}
	Recall the expressions for $f_3$ and $f_3^*$ given by \eqref{AbstractP2-f_3} and \eqref{AbstractP2-f_3*}, respectively.
	Inequality \eqref{P42P3} implies $a\le0$.
	If $x=0$, then $f_3(0)=1$, and
	inequality \eqref{P42P3} is equivalent to $f_3^*(\phi_{a})\le-1$.
	However, in view of \eqref{AbstractP2-f_3*},
	$f_3^*(\phi_{a})\ge-\frac{1}{2}$ for all $a\in\R$.
	Hence, $\partial f_3(0)=\es$.
	
	Let $0<|x|<\frac{1}{2}$.
	Then $f_3(x)=1-2|x|$.
	If $-2\le a\le0$, then $f_3^*(\phi_{a})=\frac{a}{4}$ while $ax^2-f_3(x)= ax^2+2|x|-1<\frac{a}{4}$, and consequently, inequality \eqref{P42P3} cannot be satisfied.
	If $a<-2$, then 
	$f_3^*(\phi_{a})=-1-\frac{1}{a}$, and \eqref{P42P3} becomes $ax^2+2|x|+\frac{1}{a}\ge0$, or equivalently, $(a|x|+1)^2\le0$, i.e. $a=-\frac{1}{\abs{x}}$ $(<-2)$.
	Hence, $\partial f_3(x) =\left\{\phi_{-1/{\abs{x}}}\right\}$.
	
	Let $|x|\ge\frac{1}{2}$.
	Then $f_3(x)=0$.
	If $a<-2$, then
	$f_3^*(\phi_{a})=-1-\frac{1}{a}>-\frac{1}{2}$ while $ax^2\le\frac{a}{4}<-\frac{1}{2}$, and consequently, inequality \eqref{P42P3} cannot be satisfied.
	If $-2\le a\le0$, inequality \eqref{P42P3} becomes  $\frac{a}{4}\le ax^2$, which is satisfied when either $a=0$ or $|x|=\frac{1}{2}$.
	Hence, $\partial f_3(\pm\frac{1}{2})=\left\{\phi_{a}: -2\le a\le0\right\}$, and
	$\partial f_3(x)=\{0\}$ for all $x$ with $|x|>\frac{1}{2}$.

	\section{Open questions}\label{sec-oq}
	
	\begin{enumerate}
		\item[1.] 
		In Section~\ref{AbstractExample_AC}, 
		the sum rule \eqref{AbstractF2} holds for an abstract linear space defined on a family of one dimensional functions from $\R$ to $\R$. 
		It is {interesting to investigate what the results are} for the class of functions on $\R^n$ with 
		$$
		\LL:=\{\phi_A(x): \phi_A(x):=x^TAx\}, \quad \HH:=\{\Psi_{A,b}(x)=x^TAx+b\},
		$$
		with $A$ an $n\times n$ matrix, $x\in \R^n$ and $b\in \R$. 
		Is the sum rule \eqref{AbstractF2} still holding for this more general case? What are the characterizations of $\HH-$convex functions, and $\HH-$convex sets for this case?
		\item[2.] 
		The spaces of abstract linear functions and abstract affine functions $\LL,\HH$, are equipped with pointwise convergence topology. 
		What are the semicontinuty properties of the $\eps-$subdifferential point-to set mapping?
		\item[3.] 
		Is it possible to extend the zero duality gap result (e.g., Theorem \ref{AbstractT3.6}) to an infinite sum of functions in the primal problem, {such as those considered in \cite{burachik2013strong,dinhVolle}}?
		\item[4.]
		We have shown that key constraint qualifications for zero duality gap are expressed in terms of the $\epsilon$-enlargement of the subdifferentials of the functions $f_i$. If we replace the optimization problem by the problem of finding a zero of a sum of ``abstract" maximally monotone operators, can we generate new constraint qualifications for this problem in terms of "abstract" enlargements? In other words, can  the results of Theorems \ref{AbstractT3.6} and \ref{AbstractT3} be expressed in terms of enlargements of maximally monotone operators $\{T_i\}_{i=1}^m$, for the problem of finding $x$ s.t. $0\in (\sum_{i=1}^m T_i)(x)$? To the authors' knowledge, this question is open even for the classical ``convex" case. 
	\end{enumerate}
	
	\section*{Acknowledgement(s)}
	The first author is supported by the Australian Research Council through grant IC180100030.
	The third author is supported by the Australian Research Council, project DP160100854, and benefited from the support of the 
	Conicyt REDES program 180032.
	\sloppy	
	
	\addcontentsline{toc}{section}{References}
	\bibliographystyle{tfnlm}
	\bibliography{abcon}

\begin{thebibliography}{10}
\expandafter\ifx\csname url\endcsname\relax
  \def\url#1{\texttt{#1}}\fi
\expandafter\ifx\csname urlprefix\endcsname\relax\def\urlprefix{URL }\fi
\expandafter\ifx\csname href\endcsname\relax
  \def\href#1#2{#2} \def\path#1{#1}\fi

\bibitem{Mar}
J.~E. Martínez-Legaz, Quasiconvex duality theory by generalized conjugation
  methods, Optimization 19 (1988) 1029--4945.

\bibitem{Bru}
R.~S. Burachik, A.~M. Rubinov, On abstract convexity and set valued analysis,
  J. Nonlinear Convex Anal. 9~(1) (2008) 105--123.

\bibitem{Kay}
J.~Jeyakumar, A.~M. Rubinov, Z.~Y. Wu, Generalized {F}enchel’s conjugation
  formulas and duality for abstract convex functions, J. Optim. Theory Appl.
  132 (2007) 441--458.

\bibitem{RubUde01}
A.~M. Rubinov, A.~Uderzo, On global optimality conditions via separation
  functions, J. Optim. Theory Appl. 109 (2001) 345--370.

\bibitem{YaoLi18}
C.~Yao, S.~Li, Vector topical function, abstract convexity and image space
  analysis, J. Optim. Theory Appl. 177 (2018) 717--742.

\bibitem{Iof001}
A.~Ioffe, Abstract convexity and non-smooth analysis, Adv. Math. Econ 3 (2001)
  45--61.

\bibitem{DutMarRub04}
J.~Dutta, J.~Mart\'inez-Legaz, A.~M. Rubinov, Monotonic analysis over cones:
  {II}, Optimization 53 (2004) 529--547.

\bibitem{DutMarRub04a}
J.~Dutta, J.~Mart\'inez-Legaz, A.~M. Rubinov, Monotonic analysis over cones: I,
  Optimization 53 (2004) 165--177.

\bibitem{RubDza02}
A.~M. Rubinov, Z.~Dzalilov, Abstract convexity of positively homogeneous
  functions, Journal of Statistics and Management Systems 5 (2002) 1--20.

\bibitem{NedOzdRub07}
A.~Nedi\'c, A.~Ozdaglar, A.~M. Rubinov, Abstract convexity for nonconvex
  optimization duality, Optimization 56 (2007) 655--674.

\bibitem{Shv}
A.~P. Shveidel, Abstract convex sets with respect to the class of general
  min-type functions, Optimization 52 (2003) 571--579.

\bibitem{DarMoh11}
M.~Daryaei, H.~Mohebi, Abstract convexity of extended real-valued {ICR}
  functions, Optimization 62 (2013) 835--855.

\bibitem{Rub01}
A.~M. Rubinov, Abstract convexity, global optimization and data classification,
  OPSEARCH 38~(3) (2001) 247--265.

\bibitem{EberMoh10}
A.~C. Eberhard, H.~Mohebi, Maximal abstract monotonicity and generalized
  {F}enchel's conjugation formulas, Set-Valued Anal 18 (2010) 79--108.

\bibitem{MohSam}
H.~Mohebi, M.~Samet, Abstract convexity of topical functions, J. Glob. Optim.
  58 (2014) 365--357.

\bibitem{DoaMoh09}
A.~R. Doagooei, H.~Mohebi, Monotonic analysis over ordered topological vector
  spaces: {IV}, J. Glob. Optim. 45 (2009) 355--369.

\bibitem{DutMarRub08}
J.~Dutta, J.~Mart\'inez-Legaz, A.~M. Rubinov, Monotonic analysis over cones:
  {III}, Convex Anal. 15 (2008) 581--592.

\bibitem{RubWu07}
A.~M. Rubinov, Z.~Y. Wu, Optimality conditions in global optimization and their
  applications, Math. Program., Ser. B 120 (2009) 101--123.

\bibitem{BukRub07}
R.~S. Burachik, A.~M. Rubinov, Abstract convexity and augmented {L}agrangians,
  SIAM J. Optim. 18~(2) (2007) 413--436.

\bibitem{SatMoh19-ref1}
A.~R. Sattarzadeh, H.~Mohebi, Characterizing approximate global minimizers of
  the difference of two abstract convex functions with applications, Filomat
  33~(8) (2019) 2431--2445.
\newblock \href {http://dx.doi.org/10.2298/fil1908431s}
  {\path{doi:10.2298/fil1908431s}}.

\bibitem{Wu2007-ref2}
Z.~Y. Wu, Sufficient global optimality conditions for weakly convex
  minimization problems, Journal of Global Optimization 39~(3) (2007) 427--440.
\newblock \href {http://dx.doi.org/10.1007/s10898-007-9147-z}
  {\path{doi:10.1007/s10898-007-9147-z}}.

\bibitem{MarcoRub-ref3}
M.~A. López, A.~M. Rubinov, V.~N.~V. deSerio, Stability of semi-infinite
  inequality systems involving min-type functions, Numerical Functional
  Analysis and Optimization 26~(1) (2005) 81--112.
\newblock \href {http://dx.doi.org/10.1081/NFA-200052006}
  {\path{doi:10.1081/NFA-200052006}}.

\bibitem{Ivan99-ref4}
I.~Singer, Duality in quasi-convex supremization and reverse convex
  infimization via abstract convex analysis,and applications to approximation,
  Optimization 45~(1-4) (1999) 255--307.
\newblock \href {http://dx.doi.org/10.1080/02331939908844436}
  {\path{doi:10.1080/02331939908844436}}.

\bibitem{RubGlover-ref5}
A.~M. Rubinov, B.~M. Glover, Increasing convex-along-rays functions with
  applications to global optimization, Journal of Optimization Theory and
  Applications 102~(3) (1999) 615--642.
\newblock \href {http://dx.doi.org/10.1023/A:1022602223919}
  {\path{doi:10.1023/A:1022602223919}}.

\bibitem{Rubinov1999-ref6}
A.~M. Rubinov, M.~Y. Andramonov, Minimizing increasing star-shaped functions
  based on abstract convexity, Journal of Global Optimization 15~(1) (1999)
  19--39.
\newblock \href {http://dx.doi.org/10.1023/A:1008344317743}
  {\path{doi:10.1023/A:1008344317743}}.

\bibitem{MikAnd-ref7}
M.~Andramonov, A survey of methods of abstract convex programming, Journal of
  Statistics and Management Systems 5~(1-3) (2002) 21--37.
\newblock \href {http://dx.doi.org/10.1080/09720510.2002.10701049}
  {\path{doi:10.1080/09720510.2002.10701049}}.

\bibitem{article-ref11}
G.~Crespi, I.~Ginchev, M.~Rocca, A.~Rubinov, Convex along lines functions and
  abstract convexity. i, Journal of Convex Analysis 14.

\bibitem{article-ref12}
A.~M. Rubinov, A.~Shveidel, Radiant and star-shaped functions, Pac. J. Optim
  3~(1) (2012) 193--212.

\bibitem{Ref13}
A.~M. Rubinov, I.~Singer, Topical and sub-topical functions, downward sets and
  abstract convexity, Optimization 50~(5-6) (2001) 307--351.
\newblock \href {http://dx.doi.org/10.1080/02331930108844567}
  {\path{doi:10.1080/02331930108844567}}.

\bibitem{Syg16}
M.~Syga, Minimax theorems for $\phi$-convex functions: sufficient and necessary
  conditions, Optimization 65~(3) (2016) 635--649.
\newblock \href {http://dx.doi.org/10.1080/02331934.2015.1062010}
  {\path{doi:10.1080/02331934.2015.1062010}}.

\bibitem{Syg18}
M.~Syga, Minimax theorems for extended real-valued abstract convex–concave
  functions, J. Optim. Theory Appl. 176 (2018) 306--318.
\newblock \href {http://dx.doi.org/10.1007/s10957-017-1210-4}
  {\path{doi:10.1007/s10957-017-1210-4}}.

\bibitem{EwaSyg20}
E.~M. Bednarczuk, M.~Syga, On {L}agrange duality for several classes of
  nonconvex optimization problems, in: H.~A. Le~Thi, H.~M. Le, T.~Pham~Dinh
  (Eds.), Optimization of Complex Systems: Theory, Models, Algorithms and
  Applications, Springer International Publishing, Cham, 2020, pp. 175--181.

\bibitem{GorTyk19}
V.~V. Gorokhovik, A.~S. Tykoun, Support points of lower semicontinuous
  functions with respect to the set of {L}ipschitz concave functions, Dokl.
  Nats. Akad. Nauk Belarusi 63~(6) (2019) 647--653.
\newblock \href {http://dx.doi.org/10.29235/1561-8323-2019-63-6-647-647-653}
  {\path{doi:10.29235/1561-8323-2019-63-6-647-647-653}}.

\bibitem{GorTyk20}
V.~V. Gorokhovik, A.~S. Tykoun, Abstract convexity of functions with respect to
  the set of {L}ipschitz (concave) functions, Proc. Steklov Inst. Math. 309~(1)
  (2020) S36--S46.
\newblock \href {http://dx.doi.org/10.1134/S0081543820040057}
  {\path{doi:10.1134/S0081543820040057}}.

\bibitem{Ref9}
M.~V. Dolgopolik, Abstract convex approximations of nonsmooth functions,
  Optimization 64~(7) (2015) 1439--1469.
\newblock \href {http://dx.doi.org/10.1080/02331934.2013.869811}
  {\path{doi:10.1080/02331934.2013.869811}}.

\bibitem{Rub}
A.~M. Rubinov, Abstract Convexity and Global Optimization, Kluwer Academic
  Publishers, Dordrecht, 2000.

\bibitem{ISing}
I.~Singer, Abstract Convex Analysis, Wiley-Interscience, New York, 2006.

\bibitem{Bor}
J.~Borwein, R.~S. Burachik, L.~Yao, Conditions for zero duality gap in convex
  progamming, J. Nonlinear Convex Anal. 15~(1) (2014) 167--190.

\bibitem{Bre}
H.~Brezis, Functional Analysis, Sobolev Spaces and Partial Differential
  Equations, Springer, USA, 2011.

\bibitem{Zal02}
C.~Z{\u{a}}linescu, Convex Analysis in General Vector Spaces, World Scientific
  Publishing Co. Inc., River Edge, NJ, 2002.
\newblock \href {http://dx.doi.org/10.1142/9789812777096}
  {\path{doi:10.1142/9789812777096}}.

\bibitem{Str19.2}
A.~S. Strekalovsky, Global optimality conditions and exact penalization, Optim.
  Lett. 13~(3) (2019) 597--615.
\newblock \href {http://dx.doi.org/10.1007/s11590-017-1214-x}
  {\path{doi:10.1007/s11590-017-1214-x}}.

\bibitem{LiNg08}
G.~Li, K.~F. Ng, On extension of {F}enchel duality and its application, SIAM J.
  Optim. 19~(3) (2008) 1489--1509.
\newblock \href {http://dx.doi.org/10.1137/080716803}
  {\path{doi:10.1137/080716803}}.

\bibitem{burachik2013strong}
R.~S. Burachik, S.~N. Majeed, Strong duality for generalized monotropic
  programming in infinite dimensions, Journal of Mathematical Analysis and
  Applications 400~(2) (2013) 541--557.

\bibitem{dinhVolle}
D.~T. Luc, M.~Volle, Duality for extended infinite monotropic optimization
  problems, Math. Program.\href {http://dx.doi.org/10.1007/s10107-020-01557-3}
  {\path{doi:10.1007/s10107-020-01557-3}}.

\end{thebibliography}
	
\end{document}